\documentclass[]{article}

\usepackage[T1]{fontenc}
\usepackage[cp1250]{inputenc}

\usepackage{amsmath}
\usepackage{amsfonts}
\usepackage{amssymb}
\usepackage{amsthm}
\usepackage{amscd}
\input xy
\xyoption{all}

\newcommand{\N}{\mathbb N}
\newcommand{\T}{\mathbb T}
\newcommand{\Zet}{\mathbb Z}

\newcommand{\Cs}{\mathcal C}

\newcommand{\As}{\mathcal A}
\newcommand{\Ms}{\mathcal M}
\newcommand{\Ns}{\mathcal N}
\newcommand{\Es}{\mathcal E}

\newcommand{\Bs}{\mathcal B}
\newcommand{\om}{\omega}

\font\black=cmbx10 \font\sblack=cmbx7 \font\ssblack=cmbx5 \font\blackital=cmmib10 \skewchar\blackital='177
\font\sblackital=cmmib7 \skewchar\sblackital='177 \font\ssblackital=cmmib5 \skewchar\ssblackital='177
\font\sanss=cmss10 \font\ssanss=cmss8 scaled 900 \font\sssanss=cmss8 scaled 600 \font\blackboard=msbm10
\font\sblackboard=msbm7 \font\ssblackboard=msbm5 \font\caligr=eusm10 \font\scaligr=eusm7 \font\sscaligr=eusm5
 \font\fraktur=eufm10 \font\sfraktur=eufm7 \font\ssfraktur=eufm5 

\usepackage{color,soul}
\font\bsymb=cmsy10 scaled\magstep2
\def\all#1{\setbox0=\hbox{\lower1.5pt\hbox{\bsymb
       \char"38}}\setbox1=\hbox{$_{#1}$} \box0\lower2pt\box1\;}
\def\exi#1{\setbox0=\hbox{\lower1.5pt\hbox{\bsymb \char"39}}
       \setbox1=\hbox{$_{#1}$} \box0\lower2pt\box1\;}

\def\tx#1{{\fam0\relax#1}}

\newfam\bifam
\textfont\bifam=\blackital \scriptfont\bifam=\sblackital \scriptscriptfont\bifam=\ssblackital

\newfam\blfam
\textfont\blfam=\black \scriptfont\blfam=\sblack \scriptscriptfont\blfam=\ssblack

\newfam\bbfam
\textfont\bbfam=\blackboard \scriptfont\bbfam=\sblackboard \scriptscriptfont\bbfam=\ssblackboard

\newfam\ssfam
\textfont\ssfam=\sanss \scriptfont\ssfam=\ssanss \scriptscriptfont\ssfam=\sssanss
\def\ss#1{{\fam\ssfam\relax#1}}

\newfam\clfam
\textfont\clfam=\caligr \scriptfont\clfam=\scaligr \scriptscriptfont\clfam=\sscaligr

\newfam\frfam
\textfont\frfam=\fraktur \scriptfont\frfam=\sfraktur \scriptscriptfont\frfam=\ssfraktur

\def\hpb#1{\setbox0=\hbox{${#1}$}
    \copy0 \kern-\wd0 \kern.2pt \box0}
\def\vpb#1{\setbox0=\hbox{${#1}$}
    \copy0 \kern-\wd0 \raise.08pt \box0}

\def\pmb#1{\setbox0\hbox{${#1}$} \copy0 \kern-\wd0 \kern.2pt \box0}
\def\pmbb#1{\setbox0\hbox{${#1}$} \copy0 \kern-\wd0
      \kern.2pt \copy0 \kern-\wd0 \kern.2pt \box0}
\def\pmbbb#1{\setbox0\hbox{${#1}$} \copy0 \kern-\wd0
      \kern.2pt \copy0 \kern-\wd0 \kern.2pt
    \copy0 \kern-\wd0 \kern.2pt \box0}
\def\pmxb#1{\setbox0\hbox{${#1}$} \copy0 \kern-\wd0
      \kern.2pt \copy0 \kern-\wd0 \kern.2pt
      \copy0 \kern-\wd0 \kern.2pt \copy0 \kern-\wd0 \kern.2pt \box0}
\def\pmxbb#1{\setbox0\hbox{${#1}$} \copy0 \kern-\wd0 \kern.2pt
      \copy0 \kern-\wd0 \kern.2pt
      \copy0 \kern-\wd0 \kern.2pt \copy0 \kern-\wd0 \kern.2pt
      \copy0 \kern-\wd0 \kern.2pt \box0}

\def\sT{{\ss T}}
\def\sV{{\ss V}}

\def\xd{\tx{d}}
\def\xi{\tx{i}}

\def\xv{\tx{v}}

\def\dt{\xd_{\ss T}}
\def\vt{\xv_{\ss T}}
\def\dTs{\xd_{\ss T}^*}
\def\dts{\sT^*}
\newcommand{\ee}{\end{equation}}
\newcommand{\ra}{\rightarrow}

\newcommand{\bea}{\begin{eqnarray}}
\newcommand{\eea}{\end{eqnarray}}
\newcommand{\beas}{\begin{eqnarray*}}
\newcommand{\eeas}{\end{eqnarray*}}
\newcommand{\Z}{\mathbb{Z}}
\newcommand{\R}{\mathbb{R}}

\newcommand{\ao}{\mathbf{1}^n}
\newcommand{\we}{\wedge}

\newcommand{\nn}{\nonumber}

\newcommand{\cross}{\ti}
\newcommand{\pa}{\partial}
\newcommand{\ti}{\times}
\newcommand{\A}{{\cal A}}
\newcommand{\cV}{{\cal V}}
\newcommand{\Li}{{\cal L}}

\newcommand{\cT}{{\cal T}}
\newcommand{\cF}{{\cal F}}

\newcommand{\Ll}{\Li}

\def\la{\langle}
\def\ran{\rangle}

\mathchardef\za="710B  %\alpha
\mathchardef\zb="710C  %\beta
\mathchardef\zg="710D  %\gamma
\mathchardef\zd="710E  %\delta
\mathchardef\zve="710F %\epsilon
\mathchardef\zz="7110  %\zeta
\mathchardef\zh="7111  %\eta
\mathchardef\zvy="7112 %\theta
\mathchardef\zi="7113  %\iota
\mathchardef\zk="7114  %\kappa
\mathchardef\zl="7115  %\lambda
\mathchardef\zm="7116  %\mu
\mathchardef\zn="7117  %\nu
\mathchardef\zx="7118  %\xi
\mathchardef\zp="7119  %\pi
\mathchardef\zr="711A  %\rho
\mathchardef\zs="711B  %\sigma
\mathchardef\zt="711C  %\tau
\mathchardef\zu="711D  %\upsilon
\mathchardef\zvf="711E %\phi
\mathchardef\zq="711F  %\chi
\mathchardef\zc="7120  %\psi
\mathchardef\zw="7121  %\omega
\mathchardef\ze="7122  %\varepsilon
\mathchardef\zy="7123  %\vartheta
\mathchardef\zf="7124  %\varomega
\mathchardef\zvr="7125 %\varrho
\mathchardef\zvs="7126 %\varsigma
\mathchardef\zf="7127  %\varphi
\mathchardef\zG="7000  %\Gamma
\mathchardef\zD="7001  %\Delta
\mathchardef\zY="7002  %\Theta
\mathchardef\zL="7003  %\Lambda
\mathchardef\zX="7004  %\Xi
\mathchardef\zP="7005  %\Pi
\mathchardef\zS="7006  %\Sigma
\mathchardef\zU="7007  %\Upsilon
\mathchardef\zF="7008  %\Phi
\mathchardef\zW="700A  %\Omega

\newcommand{\epf}{\hfill$\Box$}
\newcommand{\bepf}{\noindent\textit{Proof.-} }

\def\be#1{\begin{equation}\label{#1}}
\def\pd#1#2{\frac{\partial#1}{\partial#2}} %% partial  derivative
\def\wt{\widetilde}

\def\dgr{\text{\rm dgr}}

\newtheorem{re}{Remark}[section]
\newtheorem{theo}{Theorem}[section]
\newtheorem{prop}{Proposition}[section]
\newtheorem{lem}{Lemma}[section]
\newtheorem{cor}{Corollary}[section]
\newtheorem{ex}{Example}[section]
\newtheorem{definition}{Definition}[section]
\begin{document}
\title{Higher vector bundles\\ and multi-graded symplectic manifolds
\thanks{Research financed by the Polish
Ministry of Science and Higher Education %, by means of the budget for science 2006-2009,
under the grant No. N201 005 31/0115.}}

        \author{
        Janusz Grabowski$^1$, Miko\l aj Rotkiewicz$^2$\\
        \\
         $^1$ {\it Institute of Mathematics}\\
                {\it Polish Academy of Sciences}\\ {\tt jagrab@impan.gov.pl}\\\\
         $^2$ {\it Institute of Mathematics}\\
                {\it University of Warsaw}\\ {\tt mrotkiew@mimuw.edu.pl}
                }
\date{}
\maketitle
\begin{abstract}   A natural explicit condition is given ensuring that an action of the multiplicative
monoid of non-negative reals on a manifold $F$ comes from homotheties of a vector bundle
structure on $F$, or, equivalently, from an Euler vector field. This is used in showing
that double (or higher) vector bundles present in the literature can be equivalently
defined as manifolds with a family of commuting Euler vector fields. Higher vector
bundles can be therefore defined as manifolds admitting certain $\N^n$-grading in the
structure sheaf. Consequently, multi-graded (super)manifolds are canonically associated
with higher vector bundles that is an equivalence of categories. Of particular interest
are symplectic multi-graded manifolds which are proven to be associated with cotangent
bundles. Duality for higher vector bundles is then explained by means of the cotangent
bundles as they contain the collection of all possible duals. This gives, moreover,
higher generalizations of the known `universal Legendre transformation' $\dts
E\simeq\dts E^*$, identifying the cotangent bundles of all higher vector bundles in
duality. The symplectic multi-graded manifolds, equipped with certain homological
Hamiltonian vector fields, lead to an alternative to Roytenberg's picture
generalization of Lie bialgebroids, Courant brackets, Drinfeld doubles and can be
viewed as geometrical base for higher BRST and Batalin-Vilkovisky formalisms. This is
also a natural framework for studying $n$-fold Lie algebroids and related structures.

\bigskip\noindent
\textit{MSC 2000: 58A50, 53D05 (Primary); 53D17, 58C50, 17B62, 17B63, 18D05
(Secondary).}

\end{abstract}

\newpage
\section{Introduction}  The language of super-geometry is nowadays commonly used not only in some
models of mathematical physics (e.g. in Batalin-Vilkovisky formalism and topological
Quantum Field Theory \cite{AKSZ,Roy1}) or homological algebra but also for some
problems viewed earlier as purely geometrical, especially in Poisson geometry and the
theory of Lie algebroids. In this context it became evident that many canonical
super-manifolds are provided with an additional grading in the structure sheaf. In
particular, the problem of finding a proper analog of Drinfeld double Lie algebra for
Lie bialgebroids \cite{MX} and finding a nice description of Courant algebroids
\cite{LWX} (with BRST complex and the Weil algebra as particular examples) have been
solved in the language of such {\it graded (super)manifolds} by Th.~Voronov and
D.~Roytenberg \cite{Roy0,Roy,Vor1}.

On the other hand, in the traditional language of differential geometry, double (or higher) structures have
been introduced in the categorical sense. For example, {\it double vector bundles} have been understood as
"vector bundles in the category of vector bundles" (see \cite{Pr1}--\cite{Pr3}, \cite{BM,KU, Mac2})  and
recognized as the structures of great importance in the Lagrangian and Hamiltonian formalism of analytical
mechanics \cite{TU,GGU}. This has been extended for {\it double affine bundles} in \cite{GGU1}. Double
structures appeared also in symplectic and Poisson geometry with A.~Weinstein's and his collaborators work on
symplectic and Poisson groupoids \cite{CDW,We1,We2} followed by numerous works of his students and systematic
studies of K.~C.~H.~Mackenzie \cite{Mac1a}--\cite{Mac}, Y.~Kosmann-Schwarzbach \cite{YKS,YKS1} and others.

We had, however, the feeling that, on one hand, the standard definitions of a double
(or higher) vector bundle (cf. \cite{KU,Mac1a,Mac2}), although categorically nice, are
operatively too complicated, and, on the other hand, that standard concepts of
super-manifold, or even the concept of {\it N-manifold} as defined and used in
\cite{Sev, Vor1,Roy}, are still too general for many purposes. We therefore develop a
theory of higher vector bundles in the spirit of algebraically described compatibility
condition for a number of vector bundle structures and associate with them canonically
derived multi-graded super-manifolds.

Our starting point is the observation that a vector bundle can be characterized only
with the use of its homogeneous structure that leads to a much simpler definition of an
$n$-vector bundle (in classical terms). We prove namely that an $n$-vector bundle can
be equivalently characterized as a manifold with $n$ commuting Euler vector fields,
i.e. as a manifold with certain $\N$-gradation in the algebra of smooth functions. This
implies that an $n$-vector bundle, as canonically multi-graded, admits its natural
superized counterpart -- a multi-graded (super)manifold. Both concepts lead to a
unified and elegant description of various phenomena of differential geometry. Of
particular interest are symplectic multi-graded manifolds which are proven to be
associated with cotangent bundles. Duality for higher vector bundles can be explained
by means of these bundles as they contain the collection of all possible duals. In
fact, we have higher "Legendre transformations" identifying the cotangent bundles of
all these duals. The symplectic multi-graded manifolds, equipped with certain
homological Hamiltonian vector fields, lead to an alternative to D.~Roytenberg's
picture generalization of Lie bialgebroids, Courant brackets, Drinfeld doubles and can
be viewed as geometrical base for higher BRST and Batalin-Vilkovisky formalisms. This
is also a natural framework for studying $n$-fold Lie algebroids and related
structures. The paper is organized as follows.

We start with finding a simple characterization of those actions of the multiplicative
monoid $\R_+$ of non-negative reals on a manifold $F$ that come from homotheties of a
vector bundle structure on $F$. This allows us to identify a vector bundle structure
with its homogeneous structure (or, equivalently, its Euler vector field) that clearly
simplifies the whole theory, as direct comparison of the additive structures is much
more complicated. In particular, a compatibility of two vector bundle structures can be
described easily as the commutation of the corresponding Euler vector fields. We show
that this compatibility condition is equivalent to the concept of double vector bundle
described in categorical terms. In this language, a vector bundle morphism is shown to
be just a smooth map intertwining the homotheties and a vector subbundle -- as a
submanifold which is homothety-invariant.

The $n$-vector bundles $F$, whose structure is described in section 4, admit canonical
lifts of their Euler vector fields to the tangent and to the cotangent bundles $\sT F$
and $\dts F$, as we show in section 5. In particular, the iterated tangent and
cotangent bundles are canonical examples of higher vector bundles. The cotangent bundle
$\dts F$ is of particular interest, since it is canonically fibered not only over $F$
but also over all duals $F^*_{(k)}$ of $F$ with respect to all its vector bundle
structures $F\ra F_{[k]}$. The side bundles $F_{[k]}$ are canonically $(n-1)$-vector
bundles themselves. We prove the existence of a canonical identifications $\dts F\simeq
\dts F^*_{(k)}\simeq \dts F^*_{(l)}$ which are additionally symplectomorphisms. This
can be viewed as a generalization of the celebrated "universal Legendre transformation"
$\dts\sT M\simeq\dts\dts M$. Moreover, the set of higher vector bundles $\{
F,F^*_{(1)},\dots,F^*_{(n)}\}$ is closed (under natural identifications) with respect
to duality. This is a phenomenon observed first for double and triple vector bundles by
K.~Konieczna, P.~Urba\'nski and K.~C.~H.~Mackenzie \cite{KU,Mac0,Mac}.

In Section 6 we prove that symplectic $n$-vector bundles, i.e. $n$-vector bundles
equipped with a symplectic form which is linear (1-homogeneous) with respect to all
vector bundle structures, always take the form $\dts F$ for certain $(n-1)$-vector
bundle $F$. This, in turn, generalizes the known result saying that any vector bundle equipped
with a linear symplectic form is, in fact, $\dts M$.

The next two paragraphes are devoted to a natural superization of the previous notions.
In this way we get the concept, already implicitly present in the literature, of a
multi-graded manifold -- a super-manifold $\Ms$ with an $\N^n$-gradation in the
structure sheaf, and the concept of a multi-graded symplectic manifold. The crucial
here is the equivalence of categories: we have a precise prescription of passing from a
(symplectic) $n$-vector bundle to the corresponding $n$-graded (symplectic) manifold and back.

On multi-graded symplectic manifolds one can consider Master Equations, i.e. equations
of the form $\{ H,H\}=0$ for Hamiltonians of parity different from the parity of the
symplectic Poisson bracket $\{\cdot,\cdot\}$. This leads to higher multi-graded analogs
of {\it Courant algebroid} \cite{LWX,Roy} in the spirit of D.~Roytenberg's explanation
of what a Courant algebroid is. This gives also a possibility of developing higher BRST
and Batalin-Vilkovisky formalisms. The language of multi-graded manifolds is also
useful in describing the structures of $n$-fold Lie algebroids, as has been already
observed by T.~Voronov \cite{Vor2}. Section 9 is devoted to these questions together
with the concept of {\it Drinfeld $n$-tuple} -- which generalizes the notion of
Drinfeld double Lie algebra and double Lie algebroid. We end up with some results on
Drinfeld $n$-tuples, their relations to $n$-fold Lie algebroids, and examples.

To limit the size of these notes, the questions concerning the higher Dirac structures,
higher generalized geometries, etc., we postpone to a separate paper. The authors wish
to thank F.~Przytycki for helpful discussions on dynamical systems.

\section{Vector bundles and homogeneous structures}
It is a standard student exercise to show that the additive structure in a real topological vector space
determines the homogeneous structure -- the multiplication by reals. The converse is also true. The Euler's
Homogeneous Function Theorem implies that any differentiable 1-homogeneous function on $\R^n$ is automatically
linear. This suggests that the homogeneity, being much simpler notion, can be used instead linearity in
differential geometry. Let us remark that all geometric objects in this paper, like manifolds, fibrations,
etc., are assumed to be finite-dimensional, paracompact and smooth.

In this section, we will use this idea to develop a concept of a vector bundle in terms of its homogeneous
structure. To explain how we will understood the latter, let us consider a vector bundle $\zp:E\ra M$. The
homotheties in $E$ define a smooth action of the commutative monoid $(\R_+,\cdot)$ of non-negative reals,
$\R_+=\{ a\in \R: a\ge 0\}$, with multiplication:
$$h:\R_+\ti E\ra E\,,\quad h_t(e):=h(t,e)=t\cdot e\,.$$
It should be made clear that by smoothness on $\R_+$ we mean that the map can be extended to a smooth map on a
neighborhood of $\R_+$ in $\R$, thus the whole $\R$. In fact, the above $\R_+$-action can be extended to a
smooth action $\wt{h}:\R\ti E\ra E$ of the multiplicative monoid $\R$ by homotheties with possible negative
factors.

Of course, with any smooth action  $h:\R_+\ti F\ra F$, $h_t\circ h_s=h_{ts}$, of the multiplicative monoid
$(\R,\cdot)$ on a smooth manifold $F$, one can associate a smooth projection $h_0:F\ra F$ (as $h_0^2=h_0$)
onto a closed subset $N=h_0(F)$ of $F$. In this generality we can define also the {\it vertical lift}
$\cV_h:F\ra\sT F_{\mid N}$, where $\cV_h(x)\in \sT_{h_0(x)}F$ is the tangent vector at $t=0$ represented by
the smooth curve $\R_+\ni t\mapsto x_h(t):=h(t,x)\in F$. In other words,
\be{main}\cV_h:F\ra\sT F\,,\quad \cV_h(x)=\dot{x}_h(0)=\sT x_h (0,\pa_t).\ee
One can easily seen that $\cV_h(x)=0_x$ for $x\in N$. For the action by homotheties on a vector bundle we have
also the converse: $\cV_h(e_m)=0\Rightarrow e_m=0_m$.

For the terminology, note only that by a {\it vector subbundle} we always mean a subbundle over a submanifold. An important example is the vertical subbundle $\sV F_{\mid 0_M}$ in the tangent bundle $\sT F$ of a vector bundle $F$ over $M$ over the zero-section $0_M$ of $F$ which is canonically isomorphic to $F$, as shows the following.
\begin{prop}\label{p1}
For a vector bundle $F$, the vertical lift gives a canonical isomorphism of vector
bundles $\ \cV_h:F\ra\sV F_{\mid 0_M}\subset\sT F$ .
\end{prop}
\bepf In local coordinates $(x^a,y_i)$ in $F$, where $(x^a)$ are local coordinates in
$M$ and $(y_i)$ are linear coordinates in the typical fiber, we have $h(t,x,y)=(x,ty)$. In the adapted
coordinates $(x^a,y_i,\dot{x}^b,\dot{y}_j)$ in $\sT F$, the vertical lift reads $\cV_h(x,y)=(x,0,0,y)$.

\epf

\noindent Note that $\sV F_{\mid N}$ can be defined for any manifold $F$ equipped with a smooth projection
onto a subset $N$ as the subset of $\sT F_{\mid N}$ consisting of vectors which are vertical with respect to
the projection. Of course, in such generality $\sV F_{\mid N}$ need not be a vector subbundle in $\sT F$.
\begin{definition}{\rm
A {\it homogeneous structure} on a smooth manifold $F$ will be understood as a smooth action $h:\R_+\ti F\ra
F$ of the multiplicative monoid $(\R_+,\cdot)$ on $F$ which is {\it non-singular} in the sense that the
vertical lift $\cV_h(x)$ vanishes only for points $x\in N=h_0(F)$, i.e. the curves $x_h(t)$ are non-singular
for $x\notin N$.}
\end{definition}
\noindent The following theorem shows that the above property of an $\R_+$-action on $F$ determines that this
action comes from actual homotheties.
\begin{theo}\label{t1} If $h:\R_+\ti F\ra F$ is a homogeneous structure on the manifold $F$,
then there is a unique vector bundle structure on $F$ whose homotheties coincide with $h$.
\end{theo}
\bepf Working separately in components, we can assume that $F$, thus $N=h_0(F)$, is connected.
The non-singularity of $\cV=\cV_h$ (having fixed $h$ we will skip the subscript)  means that $N$ is exactly
the inverse-image by $\cV$ of the zero-section: $N=\cV^{-1}(0_F)$. The fundamental property of the vertical
lift is that it intertwines the $\R_+$-action on $F$ with the actual homotheties in $\sT F$: \be{intertwining}
\cV(h_s(x))=s\cdot\cV(x).
\ee Indeed, we get (\ref{intertwining}) from the action identity $h_t(h_s(x))=h_{ts}(x)$ after differentiating both
sides with respect to $t$ at $t=0$.

The monoid action $h$ induces a monoid representation in the tangent spaces $\sT_xF$ with  $x\in N$. To see
this, for $x\in N$, put $H_t(x):\sT_xF\ra\sT_xF$ to be the derivative $H_t(x)=D_xh_t$. It is easy to see that
$\R_+\ni t\mapsto H_t(x)$ is a representation of the monoid $(\R_+,\cdot)$ in $\sT_xF$. Indeed,
differentiating the identity $h_t\circ h_s=h_{st}$ we get $D_xh_t\circ D_xh_s=D_xh_{st}$, i.e.
\be{01}H_t(x)\circ H_s(x)=H_{ts}(x).\ee Now, put $P(x)=\frac{\xd}{\xd t}_{\mid t=0}H_t(x)$. Differentiating
(\ref{01}) with respect to $t$ at $t=0$, we get that the linear map $P(x):\sT_xF\ra\sT_xF$ commutes with
$H_s(x)$ and \be{02}P(x)\circ H_s(x)=H_s(x)\circ P(x)=s\cdot P(x).\ee Moreover, after differentiating the
latter with respect to $s$ at $s=0$, we get $P(x)^2=P(x)$. This means that $P(x)$ is a projection and that
$H_s(x)$ respects the decomposition $\sT_xF=K_x\oplus E_x$ of $\sT_xF$ into the direct sum of the kernel $K_x$
and the image $E_x$ of $P(x)$.

Let us observe that one can interpret $P(x)$ also as the vertical part of the derivative
$D_x\cV:\sT_xF\ra\sT_{\cV(x)}\sT F$ with respect to the decomposition of the space tangent to the tangent
bundle $\sT F$ at the point $\cV(x)=0_x$ of the zero-section into the vertical subspace tangent to the fiber
and the horizontal subspace tangent to the zero-section: \be{04} D_x\cV:\sT_xF\ra\sT_{\cV(x)}\sT
F=\sT_x^vF\oplus\sT_x^hF.\ee Indeed, if we trivialize locally the tangent bundle $\sT F$ in a neighborhood $U$
of $x_0\in N$ in $F$, say $\sT U=U\ti V$, $V=\sT_{x_0}F$, with coordinates $(x^a,\dot{x}^b)$, then
$H_t(x_0):V\ra V$ reads $H_t(x_0)=\pd{h_t}{x}(x_0)$ and \be{05}P(x_0)=\frac{\pa^2 h}{\pa t\pa x}(0,x_0).\ee On
the other hand, $\cV(x)=(h(0,x),\pd{h}{t}(0,x))$, so that the projection $\wt{\cV}:U\ra V$ of $\cV$ on $V$ has
the derivative
\be{06}D_{x_0}\wt{\cV}=D_{x_0}^v\cV=\frac{\pa^2 h}{\pa x\pa t}(0,x_0)=P(x_0).\ee The
family of vector space projections $P_x=P(h_0(x)):V\ra V,\quad x\in U$, in finite-dimensional vector space $V$
is locally of constant rank. Indeed, the rank of $P_x$ is the trace of $P_x$ which takes integer values and
continuously depends on $x$, thus it is locally constant. In our situation it means that the rank of the
projections $P_x$ is constant, say $k$, on $N$. By $V^1_x$ denote the image $P_x(V)$. With our local
identification, $V^1_x=E_{h_0(x)}$. The intertwining property (\ref{intertwining}) implies that $\wt{\cV}(x)$
lies in $V^1_x$. Indeed, differentiating (\ref{intertwining}) with respect to $s$ at $s=0$, we get
$D_{h_0(x)}^v(\cV(x))=P(x)(\cV(x))=\cV(x)$, i.e. \be{2} \cV(x)\in E_{h_0(x)}.\ee Since $U\ni x\mapsto P_x\in
\text{gl}(V)$ is smooth, it is clear that $P_0:=P_{x_0}$ maps $V^1_x$ isomorphically onto $V^1_0:=V^1_{x_0}$
for $x$ sufficiently close to $x_0$, say from $U$. This gives a smooth trivialization of the vector bundle
$V^1_U=\bigcup_{x\in U}V^1_x$,
$$\zF_U:V^1_U\ra U\ti V^1_0,\quad \zF(x,v_x)=(x,P_0(v_x))\,,$$
and a smooth map
$$\Psi_U:U\ra V^1_0,\quad \Psi_U=P_0\circ\wt{\cV}.$$
It is easy to see that $N\cap U=\Psi_U^{-1}\{ 0\}$ and that $\Psi_U$ is of maximal rank at points of $N$ as
the derivative $D_x\Psi_U=P_0\circ P_x$ is `onto'.

Hence, due to the Implicit Function Theorem, $N\cap U$ is a submanifold in $U$, thus the whole $N$ is a
submanifold in $F$. This implies in turn that $E=\bigcup_{x\in N}E_x$, locally isomorphic with $(N\cap U)\ti
V^1_0$, is a smooth vector subbundle in $\sT F$ over $N$. Moreover, $\cV:F\ra E$ is of maximal rank, thus a
local diffeomorphism along $N$. For, observe that any vector $v\in T_xF$ with $x\in N$, which is annihilated
by the derivative $D_x\cV$ must be annihilated by $D_x\wt{\cV}$, thus be tangent to $N$. But $\cV$ embeds $N$
as the zero-section $0_N$, so $D_x\cV$ is an injection on $T_xN\subset T_xF$. Since $\cV_{\mid N}$ is an
embedding, we can even say that $\cV$ is a global diffeomorphism on a neighborhood $U_N$ of $N$ in $F$ onto a
neighborhood $W_0$ of the zero-section in $E$. Hence, $x\mapsto s^{-1}\cdot\cV(h_s(x))$ is a diffeomorphism of
$h_s^{-1}(U_N)$ onto $s^{-1}W_0$. But, according to (\ref{intertwining}), the latter map coincides with $\cV$
which is therefore a diffeomorphism of $F=\bigcup_{s\ge 1}h_s(V_N)$ onto $E=\bigcup_{s\ge 1}sV_0$,
intertwining $h_s$ with the homothety by $s$. The vector bundle structure on $F$ can be now taken as the
pull-back of the vector bundle structure in $E$ by this diffeomorphism.

Uniqueness follows from the fact that homogeneous structure (homotheties) on a vector space completely
determines the linear structure, as 1-homogeneous smooth functions, i.e. functions satisfying $f(s\cdot
x)={s\cdot{f(x)}}$, are exactly linear functions. \epf
\begin{re} {\rm The monoid $(\R_+,\cdot)$ contains an open-dense subset of invertible
elements $(\R^*_+,\cdot)$ -- the multiplicative group of positive reals. It is clear that any action $h$ of
$(\R+,\cdot)$ restricts to a group action of $(\R^*_+,\cdot)$ which has an infinitesimal generator -- the
Euler (Liouville) vector field $\zD_h$, where $\zD_h(x)$ is the vector tangent to the curve $x_h(t)$ at $t=1$.
In the case of a homogeneous structure this is exactly the Euler (Liouville) vector field $\zD_E$ of the
vector bundle $E$. Of course, this vector field is complete and its global flow $Exp(t\zD_E)$ determines the
homogeneous structure: $Exp(t\zD_E)(x)=e^t\cdot x$. The above theorem can be reformulated in terms of this
vector field as follows. Note only that the linear part of a vector field $\zD$ on $F$ at its singular point
(zero) $x_0$ is a well-defined liner map $\sT_{x_0}F\ra\sT_{x_0}F$ which in local coordinates is represented
by the Jacobian matrix of partial derivatives of coordinates of $\zD$ near $x_0$.}
\end{re}
\begin{theo}\label{t1a}
A vector field $\zD$ on a smooth manifold $F$ is the Euler vector field of a vector bundle structure on $F$ if
and only if
\begin{description}
\item{(a)} $\zD$ is complete and the corresponding flow $\R\mapsto \zf_t=Exp(t\zD)$ of diffeomorphisms has the
limit $h_0(x)=\lim_{t\to-\infty}\zf_t$ which is a projection of $F$ onto the set $N$ of singular points of
$\zD$; \item{(b)} For every $x_0\in N$, the linear part of $\zD$ at $x_0$ is a projection.
\end{description}
\end{theo}
\bepf One can follow the idea of the above proof for $\R_+$-action but we will sketch an alternative
proof in terms of normal hyperbolicity of flows and linearization of vector fields. Since the linear part of
the vector field at singular points has only eigenvalues $0,1$, according to Shoshitaishvili Theorem, at
singular points $x_0$ of $\zD$ we have a local decompositon of the manifold into the center manifold
$W^0(x_0)$ and the unstable manifold $W^+(x_0)$. The manifold $W^0(x_0)$ is invariant, so in our case it is
unique, as it has to coincide locally with $N$. This proves that $N$ is a submanifold and we have, at least
locally, a fibration of $F$ into unstable submanifolds over $N$. But on each $W^+(x_0)$ the linear part of
$\zD$ at $x_0$ is identity, so there are no resonances and $\zD$ is smoothly equivalent to its linear part,
i.e. to the Euler vector field on $\sT_{x_0}W^0(x_0)$. These linearizations on fibers of the fibrations can be
glued together to a linearization of $\zD$ near $N$, so to a local $\R_+$-action near $N$. We can pass to the
global action thanks to the assumption (a).
 \epf

\begin{re} {\rm Of course, there are singular $(\R_+,\cdot)$-actions which therefore do not correspond to
vector bundle structures. Take for example $F=\R$ with the action $h:\R_+\ti\R\ra \R,\ (t,x)\mapsto t^2\cdot x$. It
is clear that $\cV_h$ is trivial: $\cV(x)=(0,0)\in\sT\R$ for all $x\in\R$.}
\end{re}
Theorem \ref{t1} easily implies the following.
\begin{theo}\label{t2a} Every submanifold of a vector bundle, which is invariant with
respect to homotheties, is a vector subbundle (over a submanifold of the base).
\end{theo}
\bepf Let $E$ be a submanifold of a vector bundle $F$ over $M$ which is homothety-invariant.
It is easy to see that the $\R_+$-action $h$ by homotheties, reduced to $E$, is a homogeneous structure on
$E$. This is because, clearly, $\cV_{h_{\mid E}}=(\cV_h)_{\mid E}$, since the vector tangent to a curve in a
submanifold can be naturally viewed as the vector tangent to this curve in the total manifold. This implies
that $h_{\mid E}$ is an action by homotheties with respect to a unique vector bundle structure on $E$ over the
closed submanifold $N=h_0(E)\subset M$. This vector bundle structure is a vector subbundle of $\sV E_{\mid
0_E}\subset\sV F_{\mid 0_F}\simeq F$, thus canonically a subbundle of $F$.
 \epf
\begin{re} A slightly weaker result has been communicated to us by P.~Urba\'nski who
assumed that the intersection of $E$ with every fiber of $F$ is a vector subspace.
\end{re}
It should be not surprising that the concept of a morphism in the category of vector bundles can be completely
described in terms of the corresponding  homogeneous structures.
\begin{theo}\label{t2} A smooth map $\zf:F^1\ra F^2$ between the total spaces of two vector
bundle structures $h^i_0:F^i\ra M^i$, $i=1,2$, is a morphism of the vector bundles if and only if it commutes
with homotheties
\be{mor}\zf\circ h^1_t=h^2_t\circ\zf\,.\ee
\end{theo}
\bepf Note first, that (\ref{mor}) easily implies that $\zf$ maps $M^1=h^1_0(F^1)$ into
$M^2=h^2_0(F^2)$ and fibers into fibers. We therefore can assume then that $F^i$, $i=1,2$, are just vector
spaces.

Differentiating (\ref{mor}) with respect to $t$ at $t=0$, we get
$$D_0\zf\circ\cV^1=\cV^2\circ\zf.$$
Since $\cV^i:F^i\ra\sT_0F^i$ are linear isomorphisms,
$$\zf=(\cV^2)^{-1}\circ D_0\zf\circ\cV^1:F^1\ra F^2$$ is linear.
The converse, i.e. that a vector bundle morphism commutes with homotheties is obvious.
 \epf
\begin{cor}\label{c2}
A smooth map $\zf:F^1\ra F^2$ between the total spaces of two vector bundle structures is a morphism of the
vector bundles if and only if it relates the Euler vector fields:
\be{mor1}D_x\zf(\zD_{F^1}(x))=\zD_{F^2}(\zf(x)).\ee
\end{cor}
\bepf Differentiating (\ref{mor}) with respect to $t$ at $t=1$, we get (\ref{mor1}).
Conversely, (\ref{mor1}) implies that $\zf$ intertwines the flows induced by $\zD_{F^1}$ and $\zD_{F^2}$, i.e.
$\zf\circ h^1_t=h^2_t\circ\zf$ for $t>0$, thus for all $t\in\R_+$ by continuity.

\epf

\section{Commuting Euler vector fields and double vector bundles}
Consider now two commuting homogeneous structures $h^1,h^2:\R_+\ti F\ra F$, $h^1_t\circ h^2_s=h^2_s\circ
h^1_t$ for all $s,t\in\R_+$ (or, equivalently, two commuting Euler vector fields, $[\zD^1,\zD^2]=0$). Let us
denote the corresponding bases with $E^i=h^i_0(F)$, $i=1,2$. We have, in particular,
\be{1.0} h^1_t(h^2_0(x))=h^2_0(h^1_t(x))\ee
which implies that $E^2=h^2_0(F)$ is invariant with respect to $h^1$. According to Theorem \ref{t2a}, this means
that $E^2$ is a vector subbundle of $h_0^1:F\ra E^1$ over the submanifold
\be{1.1}M=h^1_0(E^2)=h^1_0\circ h^2_0(F)=h^2_0\circ
h^1_0(F)=h^2_0(E^1)=E^1\cap E^2.\ee Analogously, $E^1$ is a vector subbundle of $h_0^2:F\ra E^2$ over $M$. We
will call them {\it side bundles}. Thus we get the following diagram of vector bundle projections
\be{db}\xymatrix{ F\ar[rr]^{h^2_0}
\ar[dd]_{h^1_0} && E^2\ar[dd]^{h^1_0} \\ \\ E^1\ar[rr]^{h^2_0} && M }\ee where we write simply $h^1_0$ also
for its restriction to $E^2$, etc. Note also that $E^1, E^2, M$ are canonically closed submanifolds in $F$ as
the zero-sections of the corresponding vector bundle structures. Moreover, according to Theorem \ref{t2},
the vertical and the horizontal arrows describe morphisms of vector bundles.

Let $F^i_x$ be the $x$-fiber in $F$ of the projection $h^i_0$. For $x\in M$, let $C_x$ be the kernel of the
linear map $h^2_0:F^1_x\ra E^2_x$. This means that $C_x$ is also the kernel of $h^1_0:F^2_x\ra E^1_x$ and
$C_x=F^1_x\cap F^2_x$. The submanifold $C_x$ of $F$ carries therefore two structures of a vector space
hereditary from $F^1_x$ and $F^2_x$ which, however, coincide according to the following proposition.
\begin{prop}\label{p2}
Two real vector space structures on a manifold $V$ with commuting homotheties coincide.
\end{prop}
\bepf Commutation of the homotheties implies that the vector space structures share the same zero 0.
Differentiating the commutation relation $h^1_t(h^2_s(x))=h^2_s(h^1_t(x))$, with respect to $t$ and $s$ at
$t=0$ and $s=0$, we get
$$D_0\cV^1(\cV^2(x))=D_0\cV^2(\cV^1(x)),$$
where $\cV^i=\cV_{h^i}$. But, for a vector space structure, $D_0\cV$ is identity on $\sT_0V$, so
$\cV^1=\cV^2$. This in turn implies $h^1=h^2$ as the vector space structure comes from $\sT_0V$ by the identification $\cV:V\ra\sT_0V$.

\epf

Let us go back to the commutative diagram of vector bundle morphisms (\ref{db}). We can reduce the whole
picture by fixing $x_0\in M$ and considering the pull-backs of $\{ x_0\}$ with respect to all the projections.
This means that we consider the situation when $M$ is just one point which can be then identified as 0 -- the
only point of the intersection of the vector spaces $E^1$ and $E^2$ as submanifolds of $F$. We know already
that $C=F^1_0\cap F^2_0$ is a common vector subspace of $F^1_0$ and $F^2_0$. We will call $C$ the {\it core}
of $(h^1,h^2)$. Since $h^1_0(E^2)=\{ 0\}$, the vector space $E^2$ is a subset, thus vector subspace, of
$F^1_0$. Analogously, $E^1$ is a subspace of $F^2_0$. Since $h^2_0$ maps the fiber $F^1_0$ linearly onto
$E^2$, its kernel $C$ is a subspace complementary to $E^2\subset F^1_0$, as $h^2_0$ is identical on $E^2$.
Thus $F^1_0=E^2\oplus C$ and, analogously, $F^2_0=E^1\oplus C$. Using trivializations of the vector bundles in
question (which always exist when the bases are contractible), we get (\ref{db}), with $M=\{ 0\}$, in the form
\be{db4}\xymatrix{
E^1\ti E^2\ti C\ar[rr]^{h^2_0} \ar[d]_{h^1_0} && E^2\ar[d]^{h^1_0} \\
E^1\ar[rr]^{h^2_0} && \{ 0\} }\ee with obvious projections which are linear maps. Note however that the
identification $F=E^1\ti E^2\ti C$ is not canonical and depends on the choice of the trivializations. Indeed,
if $(\zx_i,\zvf_a,\zvy_r)$ are linear coordinates in $E^1\ti E^2\ti C$, then a change of the bases in
$E^1,E^2$ results in a change of coordinates,
\be{cc0}(\zx'_i,\zvf'_a)=\left(\sum_j\zt_i^j\zx_j,\sum_b\zr_a^b\zvf_b\right).\ee
Further, a change of the trivialization of $h^1_0$ over $E^1$ which respects the projection on $E^2$ results
in a change of coordinates,
\be{cc1}(\zx'_i,\zvf'_a,\zvy'_r)=\left(\sum_j\zt_i^j\zx_j,\sum_b\zr_a^b\zvf_b,
\sum_b\za_r^b(\zx)\zvf_b+\sum_s\zb_r^s(\zx)\zvy_s\right).\ee For the other projection we have
\be{cc2}(\zx'_i,\zvf'_a,\zvy'_r)=\left(\sum_j\zt_i^j\zx_j,\sum_b\zr_a^b\zvf_b,
\sum_j\zg_r^j(\zvf)\zx_j+\sum_s\zd_r^s(\zvf)\zvy_s\right).\ee The changes of coordinates (\ref{cc0}) and
(\ref{cc1}) coincide if and only if they have the common form
\be{cc}\left(\zx'_i,\zvf'_a,\zvy'_r\right)=\left(\sum_j\zt_i^j\zx_j,\sum_b\zr_a^b\zvf_b,
\sum_{b,j}A_r^{bj}\zvf_b\zx_j+\sum_sB_r^s\zvy_s\right)\ee with $\zt_i^j$, $\zr_a^b$, $A_r^{bj}$, and $B_r^s$
constant. But this change of coordinates, reduced to $C$, is not linear but affine which shows that the bundle
$\zz=(h^1_0,h^2_0):F\ra E^1\ti E^2$ is affine, modelled on the trivial bundle $E^1\ti E^2\ti C$.

Let us go back to the whole generality. The collection of all $C_x$ with $x\in M$ defines a vector bundle $C$
over $M$ -- the {\it core} of $(h^1,h^2)$. If we take $x_0\in M$ and a local chart $U\subset M$ near $x_0$,
then, using local trivializations of all vector bundles over the pull-backs of $U$ (which are contractible
bases), we get from (\ref{db4}) the following local form of (\ref{db})
\be{db5}\xymatrix{
U\ti E^1_{x_0}\ti E^2_{x_0}\ti C_{x_0}\ar[r] \ar[d] & U\ti E^1_{x_0}\ar[d] \\  U\ti E^1_{x_0}\ti
E^2_{x_0}\ar[r] & U}\ee with obvious projections. One important remark is that, again, the decomposition
$$(h^2_0\circ h^1_0)^{-1}(U)=(h^1_0\circ h^2_0)^{-1}(U)=U\ti E^1_{x_0}\ti
E^2_{x_0}\ti C_{x_0}$$ depends on the choice of the trivializations. A change in trivializations results in a
change of local linear coordinates like in (\ref{cc}) but with coefficients depending on $x\in U$:
\be{ccl}\left(x'_u,\zx'_i,\zvf'_a,\zvy'_r\right)=\left(x_u,\sum_j\zt_i^j(x)\zx_j,\sum_b\zr_a^b(x)\zvf_b,
\sum_{b,j}A_r^{bj}(x)\zvf_b\zx_j+\sum_sB_r^s(x)\zvy_s\right).\ee In the coordinates
$(x_u,\zx_i,\zvf_a,\zvy_r)\in U\ti E^1_{x_0}\ti E^2_{x_0}\ti C_{x_0}$ the Euler vector fields corresponding to
the vector bundle structures $h^1$ and $h^2$ read
\be{E} \zD^1=\sum_a\zvf_a\pa_{\zvf_a}+\sum_r\zvy_r\pa_{\zvy_r}\,,\quad
\zD^2=\sum_k\zx_k\pa_{\zx_k}+\sum_r\zvy_r\pa_{\zvy_r}.\ee They clearly commute. The spaces $E^1_{x_0},
E^2_{x_0}, C_{x_0}$ can be described in terms of the Euler vector fields as submanifolds defined by equations
$\zD^1=0$, $\zD^2=0$, $\zD^1=\zD^2$, respectively. Note also that the coordinate functions
$(x_u,\zx_i,\zvf_a,\zvy_r)$ are $(\zD^1,\zD^2)$-homogeneous of bi-degrees $(0,0),(1,0),(0,1),(1,1)$,
respectively. Conversely, any change of coordinates that respects this bi-degree must be of the form
(\ref{ccl}) and it preserves $\zD^1$ and $\zD^2$. What we get locally is therefore a local form of a {\it
double vector bundle} -- the notion introduced by J.~Pradines \cite{Pr1,Pr2,Pr3} and studied in
\cite{BM,KU,Vor2}. This easily implies that also globally double vector bundles and commuting homogeneous
structures are the same objects.
Summarizing our considerations, we get the following.
\begin{theo}\label{t3a} A double vector bundle can be equivalently defined as a smooth manifold
equipped with two vector bundle structures whose Euler vector fields $\zD^1,\zD^2$ commute.
\end{theo}
\begin{theo}\label{t3}
Any double vector bundle admits an atlas with charts which are invariant with respect to both homogeneous
structures and local coordinates which are $(\zD^1,\zD^2)$-homogeneous of bi-degrees
$(0,0),(1,0),(0,1),(1,1)$. Conversely, every manifold $F$ equipped with an atlas whose charts identify some
domains in $F$ with $\prod_{i_1,i_2=0,1}V{(i_1,i_2)}$, where $V{(0,0)}$ is a domain in $\R^m$, and
$V{(1,0)},V{(0,1)},V{(1,1)}$ are $\R$-vector spaces, and the changes of coordinates respect the bi-degrees
$(i_1,i_2)$, carries a canonical structure of a double vector bundle with the Euler vector fields which are
locally of the form $\zD^1=\zD_{V(0,1)}+\zD_{V(1,1)}$ and $\zD^2=\zD_{V(1,0)}+\zD_{V(1,1)}$.
\end{theo}

\section{Higher vector bundles}
A generalization of the concept of vector bundle and double vector bundle suggested by previous considerations
is now straightforward:
\begin{definition}{\rm A smooth {\it $n$-tuple vector bundle} (shortly - {\it $n$-vector bundle}) is a smooth manifold
$F$ equipped with $n$ structures of vector bundles whose corresponding Euler vector fields $\zD^i$,
$i=1,\dots,n$, pairwise commute. A morphism between $n$-vector bundles $(F,\zD^1,\dots,\zD^n)$ and
$(F',(\zD')^1,\dots,(\zD')^n)$} is a smooth map $\zf:F\ra F'$ which relates $\zD^k$ with $(\zD')^k$, i.e.
$D_x\zf(\zD^k(x))=(\zD')^k(\zf(x))$, $k=1,\dots,n$, $x\in F$.
\end{definition}
\begin{re}{\rm
With respect to the above definition, a non-trivial permutation of the Euler vector fields leads to
non-isomorphic $n$-vector bundles. Sometimes, however, it is convenient to consider {\it weak isomorphisms},
i.e. isomorphisms up to such a permutation.}
\end{re}\noindent
An inductive reasoning, completely parallel to that proving Theorem \ref{t3}, gives the following.
\begin{theo}\label{t4}
Any $n$-vector bundle admits an atlas with charts which are invariant with respect to all the homogeneous
structures and local coordinates which are $(\zD^1,\dots,\zD^n)$-homogeneous of $n$-degrees
$i=(i_1,\dots,i_n)$, $i_k=0,1$.

Conversely, every manifold $F$ equipped with an atlas whose charts identify some domains in $F$ with
$W=\prod_{i\in\{ 0,1\}^n}V(i)$, where $V{(0)}$, $0=(0,\dots,0)$, is a domain in $\R^m$, and $V{(i)}$, $i\ne
0$, are $\R$-vector spaces, and the changes of coordinates respect the $n$-degrees $i=(i_1,\dots,i_n)$, carries
a canonical structure of an $n$-vector bundle with the Euler vector fields which are locally of the form
$$\zD^k=\sum_{i'_k\ne 0}\zD_{V(i)}\,,$$
where $i'_k=(i_1,\dots,i_{k-1},0,i_{k+1},\dots,i_n)$.
\end{theo}
It is also a straightforward inductive observation that any smooth change of coordinates in $W$ respecting the
$n$-degrees $(i_1,\dots,i_n)$ of homogeneity must be of the form
\be{vor}(v')^{j}_{i}=\sum_{\sum i^a=i}\sum_{(j_1, \ldots,j_r)} T^j_{(i^1,\ldots,i^r;j_1,\ldots,j_r)}
\prod_a v_{i^a}^{j_a}\,,
\ee
where $v_{i}^{j}$, $j=1,\dots,\text{dim}(V(i))$, are linear coordinates in $V(i)$ and $T^*_*$ are smooth
functions on $V(0)$. This shows that our $n$-bundles coincide with the $n$-tuple vector bundles described by
T.~Voronov \cite{Vor2} and the triple vector bundles (for $n=3$) studied by K.~Mackenzie \cite{Mac}.

To describe closer the structure of $n$-vector bundles, let us introduce some conventions. For $i,j\in\{
0,1\}^n$, we write $\vert i\vert=\sum_ki_k$ and $i\le j$, if $i_k\le j_k$ for all $k=1,\dots,n$. Denote also
$\zd^1=(1,0,\dots,0)$, $\zd^2=(0,1,0,\dots,0)$, etc., and $p(i)=\{ i-\zd^k\in\{ 0,1\}^n:i_k=1\}$. Let us write
also $\ao=(1,\dots,1)$ and, for $k=1,\dots,n$,  $[k]=\ao-\zd^k$.

For any $i\in\{ 0,1\}^n$ the submanifold $F_i=\bigcap_{i_k=0}\{\zD^k=0\}$ is itself an $\vert i\vert$-vector
bundle with respect to the Euler vector fields $\{\zD^k:i_k=1\}$ and bases $F_{i'}$ with $i'\in p(i)$. Thus we
get a generalization of the diagram (\ref{db}), the {\it characteristic diagram} of the $n$-vector bundle $F$,
which is a commutative diagram with $2^n$ vertices $F_i$ and vector bundle morphisms $h^k_0$ from $F_i$, with
$i_k=1$, to $F_{i'_k}$. In particular, $F$ is the total space of the $n$ vector bundle structures $h^k_0:F\ra
F_{[k]}$. The intersection of their fibers over the zero-sections gives rise to a vector bundle $C$ over the
{\it final base} $M=\bigcap_kF_{[k]}$ -- the {\it core} of the $n$-vector bundle. The final base $M$ is
locally represented by $V(0)$ and fibers of $C$ are locally represented by $V(\ao)$. The local coordinates of
$n$-degrees $\le i$ form local coordinates on $F_i$. Note that we can view formally any $n$-vector bundle as
an $(n+1)$-vector bundle by adding a trivial (zero) Euler vector field. In this way, we can regard $F_i$ as an
$n$-vector bundle with trivial Euler vector fields $\zD^k$ with $i_k=0$, i.e. with the Euler vector fields
$\zD^1,\dots,\zD^n$ from $F$ but restricted to $F_i$. Then $h^k_0$, viewed as a map from $F_i$, with $i_k=1$,
onto $F_{i'_k}$, is a morphism of $n$-vector bundles. If we remove from the characteristic diagram the total
space $F$ (together with the maps from $F$), then we get a smaller diagram of $n$-vector bundle morphisms --
the {\it base} of our $n$-vector bundle which we denote by $B=\Bs(F)$. It is easy to see that the base does
not determine $F$. There is however a final object -- the {\it base product} -- denoted by $\cross B$ such
that $\Bs(\cross B)=B$ and, for any $n$-vector bundle $F$ with base $B$, there is a submersive $n$-vector
bundle morphism $\psi_F:F\ra\cross B$ which is identical on $B$. This morphisms can be viewed as "removing the
core" operation. For example, the base product for a double vector bundle (\ref{db}) is the product (or direct
sum) of the vector bundles $E^1\ti_M E^2\simeq E^1\oplus_ME^2$. In general, $\cross\Bs(F)$ can be identified
with the image of the map $(h^1_0,\dots,h^n_0):F\ra F_{[1]}\ti\dots\ti F_{[n]}$, i.e., locally,
$\psi_F:F\ra\cross\Bs(F)$ is just the projection modulo the core:
$$\cross\Bs(F)=\prod_{i\in\{ 0,1\}^n,\, i\ne\ao}V(i).$$
The coordinate changes in $\cross\Bs(F)$ are projections of the corresponding coordinate changes (\ref{vor})
for $F$. For instance, the characteristic diagram for the triple vector bundle (cf. \cite{Mac}) looks like
\begin{equation}\label{triple} \xymatrix{
               &  F_{011}\ar[dd]\ar[dl]   &                      &  F_{111} \ar[ll] \ar[dd] \ar[dl]
  \\
F_{001}\ar[dd] &            & F_{101}\ar[ll]\ar[dd] &
  \\
               &  F_{010} \ar[ld] &                 &  F_{110} \ar[ll] \ar[ld]
  \\
F_{000}        &                   &      F_{100}\ar[ll] & }
\end{equation}
whereas its base is:
\begin{equation}\label{triple-base}
\xymatrix{
               &  F_{011}\ar[dd]\ar[dl]   &                      %& \cross\Bs(F) \ar[ll] \ar[dd] \ar[dl]
  \\
F_{001}\ar[dd] &            & F_{101}\ar[ll]\ar[dd] &
  \\
               &  F_{010} \ar[ld] &                 &  F_{110} \ar[ll] \ar[ld]
  \\
F_{000}        &                   &      F_{100}\ar[ll] & }
\end{equation}
More generally, for any $i\in\{ 0,1\}^n$, one can also define $n$-manifolds $F_{<i}=\cross\Bs(F_i)$,
which locally looks like
$$F_{<i}=\prod_{j\in\{ 0,1\}^n,\, j<i}V(j).$$
Denote by $\As^i(F)$ the space of functions on $F$ with $n$-degree $i$. It is an $\As^0(F)=C^\infty(M)$-module
which is clearly locally free and finite-dimensional, so it can be viewed as the module of sections of some
vector bundle $V^i(F)$ over $M$. For instance, for a double vector bundle $F$, the module $\As^{(1,1)}(F)$ is
locally generated by products of two coordinates of degrees $(1,0)$ and $(0,1)$, and coordinates of degree
$(1,1)$, so that the symmetric tensor product $V^{(1,0)}(F)\vee_M V^{(0,1)}(F)$ is a subbundle in
$V^{(1,1)}(F)$ and we have a short exact sequence
$$0\ra V^{(1,0)}(F)\vee_M V^{(0,1)}(F)\ra V^{(1,1)}(F)\ra C\ra 0.$$ We can
consider also the graded associative and commutative algebra $\As(F)=\bigoplus_{i\in\Z^n}\As^i(F)$ of
homogeneous functions (with the convention $\As^i(F)=\{ 0\}$ if $i\notin\N^n$). Every its part
$\As^{(i)}(F)=\bigoplus_{j\le i}\As^j(F)$ is a prototype of a {\it higher module}: we have canonical
operations $\As^j(F)\ti\As^k(F)\ra\As^{j+k}(F)$ (or $V^j(F)\otimes_M V^k(F)\ra V^{j+k}(F)$) for $j+k\le i$
with obvious properties. We can do the same with respect to the total degree and to define, for $m\in\N$, the
spaces $\As^m(F)=\bigoplus_{\vert i\vert=m}\As^i(F)$ of functions of total degree $m$, and the corresponding
higher modules $\As^{(m)}(F)=\bigoplus_{\vert i\vert\le m}\As^i(F)$. They correspond to certain vector bundles
$V^m(F)$ and $V^{(m)}(F)$ over $M$.

\section{The tangent lift, the phase lift, and duality}
In this section we show how to lift Euler vector fields to the tangent and the cotangent bundle. For the
tangent and cotangent lifts of vector fields we refer to \cite{YI,GU,GU0}. Note only that both lifts respect
the Lie bracket.

Applying the tangent functor to homotheties associated with a homogeneous structure $h:\R_+\ti E\ra E$ we get a
new homogeneous structure $\dt h$, $(\dt h)_t=\sT (h_t)$. Indeed, $\dt h:\R_+\ti \sT E\ra \sT E$ is clearly a
smooth action of $(\R_+,\cdot)$ and the non-singularity assumption is preserved. In the adapted local
coordinates in $\sT E$ we have
\be{tl} \dt h(t,x,y,\dot{x},\dot{y})=(x,t y,\dot{x},t\dot{y}).\ee
Thus the projection $(\dt h)_0$ maps $\sT E$ onto $\sT M$ and the corresponding Euler vector field is the {\it
(complete) tangent lift} of the Euler vector field of $h$,
\be{tE} \dt \zD_E=\sum_ky_k\pa_{y_k}+\sum_k\dot{y}_k\pa_{\dot{y}_k}.\ee
Note that the tangent lift $\dt \zD_E$ is linear, i.e. commutes with the Euler vector field of the tangent
bundle $\zD_{\sT E}$, and on $E$ it reduces to $\zD_E$. Moreover, the tangent lifts of commuting vector fields
commute, so we get the following.
\begin{theo}\label{te1} The tangent bundle of an $n$-vector bundle $(F,\zD^1,\dots,\zD^n)$ is
canonically an $(n+1)$-vector bundle with respect to the Euler vector fields
$\dt\zD^1,\cdots,\dt\zD^n,\zD_{\sT F}$. The corresponding side bundles are $F$ and $\sT F_{[k]}$,
$k=1,\dots,n$, respectively. In particular, the iterated tangent bundle $\sT^{(n)}M=\sT\sT\cdots\sT M$ is
canonically an $n$-vector bundle with $(\sT^{(n)}M)_i\simeq(\sT^{(\vert i \vert)}M)$.
\end{theo}\noindent
If $\prod_{i\in\{ 0,1\}^n}V(i)$ are local charts in $F$ as in Theorem \ref{t4}, then we have
%%$V_{\sTF}(0,1)=\sT V(0)$ and
$V_{\sT F}((i,0))=V_{\sT F}((i,1))=V(i)$, $i>0$, for factors of local charts in $\sT F$.
In particular, $(\sT F)_{(i,0)}=F_i$ and $(\sT F)_{(i,1)}=\sT F_i$.

The "phase functor" has not as good properties as the tangent one, since, in general, it associates only
relations with smooth maps. Therefore the cotangent lift $\dTs \zD_E$ of $\zD_E$, which by definition is the
hamiltonian vector field of the linear function $\zi_{\zD_E}$ on $\sT^*E$ represented by $\zD_E$, is not an
Euler vector field. We get, however, an Euler vector field, denoted by $\dts \zD_E$ and called the {\it phase
lift} of $\zD_E$, if we add the Euler vector field $\zD_{\sT^* E}$ of the cotangent bundle,
\be{plift}\dts \zD_E=\dTs \zD_E+\zD_{\sT^* E}\,.\ee In the adapted local coordinates,
$$\dTs\zD_E(x,y,p,\zp)=\sum_ky_k\pa_{y_k}-\sum_j\zp_j\pa_{\zp_j}$$ and
\bea\label{cl} \dts \zD_E(x,y,p,\zp)&=&(\sum_ky_k\pa_{y_k}-\sum_j\zp_j\pa_{\zp_j})+(\sum_j\zp_j\pa_{\zp_j}+\sum_ap_a\pa_{p_a})\\
&=& \sum_ky_k\pa_{y_k}+\sum_ap_a\pa_{p_a}.\nn\eea The Euler lift $\dts \zD_E$ is linear, i.e. commutes with
$\zD_{\sT^* E}$, and on $E$ it reduces to $\zD_E$. The base of the corresponding homogeneous structure $\dts
h$ is canonically identified with the dual bundle $E^*$ which is canonically embedded in $\dts E$. The
homogeneous structure $\dts h$ will be called the {\it phase lift} of $h$. It commutes with the standard
homogeneous structure $h^{\dts E}$ on $\dts E$. Thus the cotangent bundle $\dts E$ is canonically a double
vector bundle with respect to the pair of commuting Euler vector fields $(\dts\zD_E,\zD_{\dts E})$. It is well
known that there is  a canonical isomorphism of double vector bundles (cf. \cite{Du, KU, GU2}) being
simultaneously a symplectomorphism of the canonical symplectic structures:
\be{alpha}\xymatrix{
 & \sT^\ast E \ar[rrr]^{\cT_E} \ar[dr]^{h^{\dts E}_0}
 \ar[ddl]_{(\dts h)_0}
 & & & \sT^\ast E^*\ar[dr]^{(\dts h^{E^*})_0}\ar[ddl]_/-20pt/{h^{E^*}_0}
 & \\
 & & E\ar[rrr]^/-20pt/{-id}\ar[ddl]_/-20pt/{h_0}
 & & & E \ar[ddl]_{h_0}\\
 E^\ast\ar[rrr]^/-20pt/{id}\ar[dr]^{h^{E^*}_0}
 & & & E^\ast\ar[dr]^{h^{E^*}_0} & &  \\
 & M\ar[rrr]^{id}& & & M &
}\ee which in local coordinates reads
\be{leg}\cT_E(x,y,p,\zp)=(x,\zp,p,-y)
\ee
and identifies $(\zD_{\dts E},\dts\zD_E)$ with $(\dts\zD_{E^*},\zD_{\dts E^*})$. This isomorphism, called
sometimes a {\it Legendre transform}, has been first discovered by W.~M.~Tulczyjew \cite{Tu} for $E=\sT M$ in
the context of Legendre transformation in analytical mechanics. Since $\cT_E$ is a symplectomorphism, we get
additionally that the canonical symplectic form $\zw_E$ on $\dts E$ is 1-homogeneous not only with respect to
$\zD_{\dts E}$ but also with respect to the phase lift $\dts\zD_E$. These properties completely determine the
vector field $\dts\zD_E$ if its restriction to $E$ is given. Namely, we have the following.
\begin{prop}\label{p3} Any vector field $X$ on $\dts M$ which commutes with the Euler vector field $\zD_{\dts
M}$ and satisfies $\Ll_X\zw_M=a\cdot\zw_M$, where $a\in\R$ and $\zw_M$ is the canonical symplectic form on
$\dts M$, is tangent to $M$ and completely determined by $a$ and its restriction to $M$. In particular, the
cotangent lift $\dTs Y$ is the unique linear and hamiltonian extension to $\dts M$ of a vector field $Y$ on
$M$.
\end{prop}
\bepf Write $X=\sum_j\left(f_j(x,p)\pa_{x^j}+g_j(x,p)\pa_{p_j}\right)$ in local Darboux coordinates
$(x^j,p_k)$. The property $[X,\zD_{\dts M}]=0$ implies easily that $f_j$ are of 0-homogeneous and $g_j$ are
1-homogeneous with respect to $\zD_{\dts M}$, i.e. $f_j(x,p)=f_j(x)$ and $g_j(x,p)=\sum_kg_j^k(x)p_k$. Now,
$$a\cdot\sum_j\xd p_j\we\xd x^j=\Ll_X\zw_M=\sum_{j,k}\left(g_j^k(x)+\frac{\pa f_k}{\pa x^j}(x)\right)\xd p_k\we\xd
x^j,$$ i.e.
$$g_j^k(x)=a\cdot\zd^k_j-\frac{\pa f_k}{\pa x^j}(x).$$
Thus, for a given $a$, the vector field $X$ is completely determined by its restriction
$\sum_jf_j(x)\pa_{x^j}$ to $M$. \epf

\medskip
Since the phase lift of an Euler vector field is Euler, the cotangent bundle of an $n$-vector bundle
$(F,\zD^1,\dots,\zD^{n})$ is canonically an $(n+1)$-vector bundle with respect to the Euler vector fields
$\dts\zD^1,\cdots,\dts\zD^{n},\zD_{\sT^* F}$ due to the following proposition.
\begin{prop}\label{p3a} The phase lifts $\dts X$ and $\dts Y$ of vector fields $X,Y$ on $M$ commute if and only if $X$
and $Y$ commute.
\end{prop}
\bepf The linear functions $\zi_X,\zi_Y$ on $\dts M$, corresponding to commuting vector fields $X$ and $Y$,
commute with respect to the symplectic Poisson bracket, so that they hamiltonian vector fields $\dTs X$ and
$\dTs Y$ commute. The cotangent lifts $\dTs X$ and $\dTs Y$ are linear vector fields on $\dts M$, so they
commute with the Euler vector field $\zD_{\dts M}$. Hence
$$[\dts X,\dts Y]=[\dTs X+\zD_{\dts M},\dts Y+\zD_{\dts M}]=0.$$ Conversely, if $\dts X$ and $\dts Y$ commute,
then $[X,Y]=[\dts X,\dts Y]_{\mid M}=0$. \epf

\medskip In homogeneous local coordinates $(x^j)$ on $F$, put $g_k(x^j)$ to be the degree of $x^j$ with respect to
$\zD^k$, $\zD^k(x^j)=g_k(x^j)x^j$. Then, in the adapted local coordinates $(x^j,p_s)$ in $\dts F$,
\be{pl}\dts\zD^k=\sum_kg_k(x^j)x^j\pa_{x^j}+\sum_k(1-g_k(x^j))p_j\pa_{p_j}.\ee
If $\prod_{i\in\{ 0,1\}^n}V(i)$ are local charts in $F$ as in Theorem \ref{t4}, then we have $V_{\dts
F}(\ao,1)=\dts V(0)$ and
$$V_{\dts F}((i,0))=V(i),\quad V_{\dts F}((i,1))=V(\ao-i)^*,\
i>0,$$ for factors of local charts in $\dts F$. In particular, $(\dts F)_{(i,0)}=F_i$. The side bundles of the
$(n+1)$-vector bundle $\dts F$ are $F$ and $F^*_{(k)}$, $k=1,\dots,n$, where
$F^*_{(k)}=(F,\zD^1,\dots,\zD^n)^*_{\zD^k}$ is the vector bundle dual to the vector bundle structure
$h^k_0:F\ra F_{[k]}$ determined by the Euler vector field $\zD^k$, respectively. But, according to
(\ref{alpha}), $\dts F\simeq\dts F^*_{(k)}$, $\dts\zD^k\simeq\zD_{\dts F^*_{(k)}}$, so up to all these identifications we can write
\be{alpha1}(\dts F,\dts\zD^1,\dots,\dts\zD^{n},\zD_{\dts F},)\simeq (\dts
F^*_{(k)},\zD_{\dts F^*_{(1)}},\dots,\zD_{\dts F^*_{(n+1)}}),
\ee
where we use the convention $F=F^*_{(n+1)}$ and $\zD_{\dts F^*_{(n+1)}}=\zD_{\dts F}$. The dual bundle
$F^*_{(k)}$ is canonically an $n$-vector bundle (with respect to the restrictions of the corresponding Euler
vector fields):
$$(F^*_{(k)},\zD_{\dts
F^*_{(1)}},\dots,\zD_{\dts F^*_{(k-1)}}, \zD_{\dts F^*_{(k+1)}},\dots,\zD_{\dts F^*_{(n+1)}}).$$ One can also
easily derive the fact that the set of $n$-vector bundles
$$(F,\zD^1,\dots,\zD^n)^*=\{ F,F^*_{(1)},\dots,F^*_{(n)}\},$$
the {\it set of duals} of the $n$-vector bundle $(F,\zD^1,\dots,\zD^n)$, is closed with respect to passing to
the dual bundle with respect to any of the vector bundle structure on them (cf. \cite{KU,Mac}). The
corresponding isomorphisms respect the $n$-bundle structures, if we accept the weak isomorphisms related to
reordering of the Euler vector fields (or, better to say, by fixing the original order $\zD_{\dts
F^*_{(1)}},\dots,\zD_{\dts F^*_{(n+1)}}$). In fact, for $k,l=1,\dots,n$, $l\ne k$,
\bea\label{duals}&(F^*_{(k)},\zD_{\dts F^*_{(1)}},\dots,\zD_{\dts F^*_{(k-1)}}, \zD_{\dts
F^*_{(k+1)}},\dots,\zD_{\dts F^*_{(n+1)}})^*_{\zD_{\dts F^*_{(l)}}} =
\left(\prod_{i_k=0}V(i)\ti\prod_{i_k=1}V(i)^*\right)^*_{\zD_{\dts F^*_{(l)}}}=\nn\\
&\left(\prod_{i_k,i_l=0}V(i)\ti\prod_{i_k=1,(\ao-i)_l=0}V(i)^*\right)\ti
\left(\prod_{i_k=0,i_l=1}V(i)^*\ti\prod_{i_k=1,(\ao-i)_l=1}V(i)\right)=\nn\\
&\prod_{i_l=0}V(i)\ti\prod_{i_l=1}V(i)^*= (F^*_{(l)},\zD_{\dts F^*_{(1)}},\dots,\zD_{\dts F^*_{(l-1)}},
\zD_{\dts F^*_{(l+1)}},\dots,\zD_{\dts F^*_{(n+1)}})\,.\eea This implies that the set of duals of an
$n$-vector bundle contains, in general, $(n+1)$-elements up to isomorphisms, if we accept weak isomorphisms, or $(n+1)!$ elements, if we count
permutations of the $n$-vector bundle structures. For example, if we start with a double vector bundle
(\ref{db}) with the core $C$, then we get the following triple vector bundle: \[ \xymatrix{
               &  F^*_{(2)} \ar[dd]\ar[dl]   &                     &  T^*F \ar[ll] \ar[dd] \ar[dl]
  \\
C^* \ar[dd] &            & F^*_{(1)} \ar[ll]\ar[dd] &
  \\
               &  E^2 \ar[ld] &                 &  F \ar[ll] \ar[ld]
  \\
M        &                   &      E^1\ar[ll] & }
\]
Our observation can be summarized as follows.
\begin{theo}
The cotangent bundle $\dts F$ of an $n$-vector bundle $(F,\zD^1,\dots,\zD^n)$ is canonically an $(n+1)$-vector
bundle with respect to the Euler vector fields $\dts\zD^1,\cdots,\dts\zD^n,\zD_{\sT^* F}$ and the side bundles
$F$ and $F^*_{(1)},\dots,F^*_{(n)}$-the dual bundles of $F$ with respect to all the vector bundle structures
on $F$. There are canonical isomorphisms of the $(n+1)$-vector bundles $\dts F\simeq\dts F^*_{(k)}$. Moreover,
the duals of the $n$-vector bundle $F^*_{(k)}$ are canonically isomorphic to
$F,F^*_{(1)},\dots,F^*_{(k-1)},F^*_{(k+1)},\dots,F^*_{(n)}$. In particular, the iterated cotangent bundles
$(\sT^*)^{(n)}M=\sT^*\sT^*\cdots\sT^* M$ are canonically $n$-vector bundles.
\end{theo}

\section{Symplectic and Poisson $n$-vector bundles}

\begin{definition}{\rm A {\it symplectic $n$-vector bundle} is an $n$-vector bundle
$(F,\zD^1,\dots,\zD^n)$ equipped with a symplectic form $\zW$ which is 1-homogeneous with respect to all
vector bundle structures:
\be{svb} \Ll_{\zD^k}\zW=\zW,\quad k=1,\dots,n\,,\ee
where $\Ll$ denotes the Lie derivative.}
\end{definition}
An example of a symplectic vector bundle is the cotangent bundle $\sT^* M$ with the canonical symplectic form
$\zw_M$.  Consequently, the cotangent bundle of any $(n-1)$-vector bundle $(E,\zD^1,\dots,\zD^{n-1})$ is a
canonical example of a symplectic $n$-vector bundle. Indeed, we know already that the canonical symplectic
structure $\zw_M$ is 1-homogeneous with respect to $\zD_{\dts F}$ and with respect to any phase lift.
\begin{theo}\label{tsym} Any symplectic $n$-vector bundle $(F,\zD^1,\dots,\zD^n,\zW)$, $n\ge 1$, is canonically
isomorphic to the cotangent bundle over each of its side bundles $F_{[k]}$, equipped with the canonical
symplectic form:
$$(F,\zD^1,\dots,\zD^n,\zW)\simeq(\dts F_{[k]},\dts(\zD^1_{\mid F_{[k]}}),\dots,
\dts(\zD^{k-1}_{\mid F_{[k]}}),\zD_{\dts F_{[k]}}, \dts(\zD^{k+1}_{\mid
F_{[k]}}),\dots,\dts(\zD^n_{\mid F_{[k]}}),\zw_{F_{[k]}}),$$ $k=1,\dots,n$. In
particular, all symplectic $n$-vector bundles
$$(\dts F_{[k]},\dts(\zD^1_{\mid F_{[k]}}),\dots,\dts(\zD^{k-1}_{\mid F_{[k]}}),\zD_{\dts F_{[k]}},
\dts(\zD^{k+1}_{\mid F_{[k]}}),\dots,\dts(\zD^n_{\mid F_{[k]}}),\zw_{F_{[k]}})$$
$k=1,\dots,n$, are canonically isomorphic.
\end{theo}
\bepf Since $\zW$ is a 1-homogeneous symplectic form on the vector bundle $h^k_0:F\ra
F_{[k]}$, we have a canonical isomorphism $\zf_k:(F,\zW)\ra(\dts
F_{[k]},\zw_{F_{[k]}})$ of symplectic vector bundles which identifies $\zD^k$ with
$\zD_{\dts F_{[k]}}$. But $F_{[k]}$ is an $(n-1)$-vector bundle with respect to the
restrictions of $\zD^1,\dots,\zD^{k-1},\zD^{k+1},\dots,\zD^n$, so ${\dts F_{[k]}}$,
thus $F$, is a symplectic $n$-vector bundle with respect to the Euler vector fields
$$\dts(\zD^1_{\mid F_{[k]}}),\dots,\dts(\zD^{k-1}_{\mid F_{[k]}}),\zD_{\dts F_{[k]}},
\dts(\zD^{k+1}_{\mid F_{[k]}}),\dots,\dts(\zD^n_{\mid F_{[k]}}).$$ Since $\zf_k$ is identity on $F_{[k]}$, the
vector field $\dts(\zD^{j}_{\mid F_{[k]}})$ coincides with $(\zf_k)_*(\zD^j)$ on $F_{[k]}$. But the linear
vector field $X=(\zf_k)_*(\zD^j)$ on the cotangent bundle $\dts F_{[k]}$, which additionally satisfies
$\Ll_X\zw_{F_{[k]}}=\zw_{F_{[k]}}$ is completely determined by its values on $F_{[k]}$, so
$$\dts(\zD^{j}_{\mid F_{[k]}})=(\zf_k)_*(\zD^j).$$
\epf
\begin{definition}{\rm A {\it Poisson $n$-vector bundle} is an $n$-vector bundle
$(F,\zD^1,\dots,\zD^n)$ equipped with a Poisson tensor $\zL$ which is linear, i.e. homogeneous of degree -1,
with respect to all vector bundle structures:
\be{pvb} \Ll_{\zD^k}\zL=-\zL,\quad k=1,\dots,n\,.\ee
}
\end{definition}
\noindent Of course, any symplectic $n$-vector bundle is automatically a Poisson $n$-vector bundle. Since any
linear Poisson structure on a vector bundle $E$ corresponds to certain de Rham derivative in the Grassmann algebra
$A(E^*)$ associated with the dual bundle, we can associate with any Poisson $n$-vector bundle $F$ the de Rham
derivatives $\xd_k$ in $A(F^*_{(k)})$. For Poisson structures, homogeneity of degree -1 we call linearity,
since the corresponding Poisson bracket is closed on linear (1-homogeneous) functions. This is exactly the Lie
algebroid bracket on section of the dual bundle. Thus we can state the following.
\begin{prop} Any Poisson $n$-vector bundle $F$ induces Lie algebroid structures on all dual vector bundles
$F^*_{(k)}\ra F_{[k]}$.
\end{prop}\noindent
The Lie algebroid structures on all duals of an $n$-vector bundle $F$ we will call {\it concordant}, if they
are obtained in the above way -- from a Poisson $n$-vector bundle structure on $F$.
\section{Multi-graded manifolds}
A graded manifold is a super-manifold equipped with an additional grading in the structure sheaf. The
coordinate transformations are required to preserve this grading. The calculus on graded manifolds has been
developed e.g. in \cite{Vor1,Roy,Meh}. Our aim is to describe super-manifolds graded by $n$-tuples of
non-negative integers (i.e. by $\N^n$).
\begin{definition}{\rm Let $G$ be an abelian semigroup, $G\ni g\mapsto p_g\in\N, g\in G$, be any function such that
$p_g\neq 0$ only for finitely many $g\in G$, and let $g\mapsto \tilde{g}$ be a semigroup homomorphism
$G\to\Z_2$. A {\it $G$-graded manifold\ } $\Ms$ of dimension $(p)$ is a super-manifold whose local coordinates
$(x_i)$ can be chosen homogeneous with respect to a $G$-gradation in the structure sheaf which agrees with the
$\Z_2$-gradation, i.e. such that the $G$-degrees coincide with the parity:
$$x_ix_j=(-1)^{\wt{g}(x^i)\wt{g}(x_j)}x_jx_i,$$ and the changes of coordinates respect the gradation.}
\end{definition}
Note also that, if $G$ is an abelian monoid (with additive notation) then the local coordinates of degree
$0\in G$ give rise naturally to a graded submanifold $M$ of $\Ms$ which is a standard (non-graded) smooth
manifold together with a projection $\Ms\ra M$. For other concepts of graded differential geometry we refer to
\cite{Roy,Vor1} or to \cite{Meh}. Especially, the concept of degree-shifting functor $[h]$ we borrow from the
latter. Let $h\in G$. The degree-shifting functor $[h]$ acts on the category of $G$-graded vector spaces and
assigns to a $G$-graded space $V=\oplus_{g\in G}V_g$ the space $W=\oplus_{g\in G}W_{g+h}$, where $W_{g+h}$
consists of the same elements as $V_g$ but has degree $g+h$. Any functor on $G$-graded vector spaces gives
rise to an operation on $G$-graded vector bundles. Note that the shift operator $[h]$ has the effect of
decreasing the degree of fiber coordinates of  a $G$-graded vector bundle $\Es\to \Ms$ by $h\in G$.
In all our cases $G$ will be the group
$\Zet^n$ (or its sub-semigroup) and $g=(g_1, \ldots,g_k) \mapsto \tilde{g} = (g_1+\ldots +g_k) \, mod\, 2$.

For a $G$-graded manifold $\Ms$ we denote by $\As(\Ms)=\bigoplus_{g\in G}\As^g(\Ms)$
the $G$-graded algebra of polynomial functions on $\Ms$.

\begin{definition}{\rm An $n${\it -graded manifold} is an $\N^n$-graded manifold $\Ms$ which
admits an atlas with local coordinates of degrees $\le \ao=(1,\dots,1)\in\N^n$. }
\end{definition}\noindent
Similarly as in the case of $n$-vector bundles, we have the algebra $\As(\Ms)=\bigoplus_{i\in\N^n}\As^i(\Ms)$
of polynomial functions. The difference is that this graded associative algebra is graded commutative instead
of being just commutative. The $\As^0(\Ms)=C^\infty(M)$-modules $\As^i(\Ms)$ and the higher modules
$\As^{(i)}(\Ms)$, as well as the corresponding vector bundles $V^i(\Ms)$ and $V^{(i)}(\Ms)$ are defined
completely analogously. We can also pass to to the corresponding objects with respect to the total degree.
\begin{re}{\rm Passing from $\N^n$- or $\Z^n$-degree $i$ to the {\it total degree} $\vert
i\vert=\sum_ki_k$ allows us to associate with every $n$-graded manifold an {\it $N$-manifold of degree $n$},
in the terminology introduced by P.~\v{S}evera \cite{Sev} and exploited by D.~Roytenberg \cite{Roy0,Roy}. }
\end{re}
A convenient way to describe the $\N^n$-gradation in an $n$-graded manifold $\Ms$ is to consider the Euler
vector fields $\zD^k_\Ms$, $k=1,\dots,n$, whose eigenvalues represent the degrees of homogeneous functions
$g(f)=(g_1(f),\dots,g_n(f))$. In local coordinates $(x^j)$,
\be{eu}\zD^k_\Ms=\sum_jg_k(x^j)x^j\pa_{x^j},\ee so
$f$ is of degree $i\in\N^n$ if $\zD^k_\Ms(f)=i_kf$, $k=1,\dots,n$. We have a fundamental correspondence
between $n$-vector bundles and $n$-graded manifolds.
\begin{theo} With every $n$-vector bundle $F=(F,\zD^1,\dots,\zD^n)$ one can canonically associate an $n$-graded
manifold $\Ms_F=\text{\rm gr}(F,\zD^1,\dots,\zD^n)$ such that local coordinates in $F$ of $n$-degree $i\le\ao$
correspond to graded local coordinates in $\Ms_F$ of degree $i$. This correspondence gives an equivalence of
the corresponding categories.
\end{theo}
\bepf
Assume that an $n$-vector bundle $F$ is given by an atlas in which local coordinates
$v_i^j$'s transform as in (\ref{vor}). Passing from $F$ to a super-manifold structure
requires a slight caution because in general the transformation formula (\ref{vor})
does not work in a super-manifold context.

Let us introduce coordinates $\theta_i^j$, of degree $i\in\{0,1\}^n$  on a domain $V_{(0)}\subset\R^m$,
corresponding to the coordinates $v_i^j$. Let us fix an order $\prec$ on the set $\{0,1\}^n$ such that
$\zd^1\prec\ldots\prec \zd^n$. Let $J(i) = (j_1,\ldots,j_{|i|})$ be the growing sequence of those
$k=1,\dots,n$, for which $i_k=1$. Let $[i^1,\ldots, i^r]\in \{\pm1\}$ be the sign of the permutation $(J^1,
\ldots, J^r)$ of the set $J(i)=J^1\cup\ldots\cup J^r$, where $J^a = J(i^a)$ and $i=i^1+\ldots+i^r$. We claim
that the following change of coordinates
\be{vorgr}
(\theta')^{j}_{i}=\sum_{\scriptstyle \sum_ai^a=i }\ \sum_{(j_1, \ldots,j_r)} [i^1,\ldots, i^r]\,
T^j_{(i^1,\ldots,i^r;j_1,\ldots,j_r)}\, \theta_{i^1}^{j_1} \ldots \theta_{i^r}^{j_r}\,,
\ee
where $\scriptstyle 0<i^1\prec\ldots\prec i^r$, satisfies the cocycle condition. This is so because of the
following easy properties of the introduced sign
\be{a}
 [i^{\sigma_1},\ldots, i^{\sigma_r}]\theta_{i^{\sigma_1}}^{j_{\sigma_1}} \ldots \theta_{i^{\sigma_r}}^{j_{\sigma_r}} =
 [i^1,\ldots, i^r] \theta_{i^1}^{j_1} \ldots \theta_{i^r}^{j_r}
\ee
for any permutation $\sigma$ and
\be{b}
 [i^{1,1}+\ldots +i^{1,s}, i^2,\ldots, i^r] [i^{1,1},\ldots, i^{1,s}] = [i^{1,1},\ldots ,i^{1,s}, i^2,\ldots, i^r].
\ee
One can obtain the corresponding $n$-graded manifold $\Ms_F$ also by means of applying the degree-shifting
functors (\cite{Meh}, Proposition 2.2.27).

Conversely, given a $n$-graded manifold one easily recovers the non-graded transformation functions $T^*_*$,
which produce the $n$-vector bundle $F$. Also the morphisms in the considered categories are in one-to-one
correspondence if we apply the rules analogous to the rules just described for the coordinate changes. In
particular, the rule (\ref{vorgr}) describes graded diffeomorphisms.\epf

\medskip\noindent
Similarly like in the $n$-vector bundle case, any $n$-graded manifold gives rise to a commutative diagram of
graded vector bundle projections $\zt_{(i,i')}:\Ms_i\ra\Ms_{i'}$, where $i\in\N^n$, $i'\in p(i)$, and $\Ms_i$
is the graded submanifold of $\Ms$ with local coordinates reduced to those whose degrees are $\le i$. Of
course, $\Ms_i$ is canonically an $\vert i\vert$-graded manifold. In this way we get graded vector bundles
$\zt^k:\Ms\ra\Ms_{[k]}$ and the dual bundles $\zs^k:\Ms^*_{(k)}\ra\Ms_{[k]}$.

Given an $n$-graded manifold $\Ms$ let us denote by ${\rm dgr}(\Ms)$ the degradation of $\Ms$, i.e. an
$n$-vector bundle $F$ such that ${\rm gr}(F) = \Ms_F$. If $A$ is a subset of $\{1, \ldots, n\}$, $\#A = k$,
we can consider $F$ as $(n-k)$-vector bundle with respect to the Euler vector fields $\zD^s$ with $s\not\in A$.
We denote it by $_AF$ and define $_A\Ms := {\rm dgr}(\Ms, \zD^A)$ as $(n-k)$-graded manifold associated with $_AF$,
i.e. ${\rm gr}(_AF, \zD^A) = {_A\Ms}$. Note that the final base of  ${_AF}$ (and so the support of ${_A\Ms}$)
has $\#A$-vector bundle structures. If $A$, $B$ are disjoint subsets of $\{1, \ldots, n\}$ then
$_{A\cup B}\Ms = {_A(_B\Ms)}$, since both sides are $(n-\#A - \#B)$-graded manifolds associated with
the $(n-\#A - \#B)$-vector bundle $_{A\cup B}F$.

\begin{ex}{\rm With a vector bundle $E$ over $M$ we associate the $\N$-graded manifold $\Ms_E$. Local
coordinates $(x^a)$ on $M$ and a basis of local sections of the dual bundle $E^*$ give rise to local
homogeneous coordinates $(x^a,y^i)$ on $E$. The local coordinates on $\Ms_E$ are $(x^a,\zx^i)$ of degrees,
respectively, 0 and 1, and the same change of coordinates as described by the definition of the vector bundle
$E^*$. Thus, with every $i$-section $\zn$ from the Grassmann algebra $A(E^*)$ we associate a function
$\zi_\zn$ of degree $i$ on $\Ms_E$ in an obvious way, so that $\As(\Ms_E)\simeq A(E^*)$. A Lie algebroid
structure on $E$ is the same as a linear Poisson structure on $E^*$, or the same as a homological vector field
$Q$, $[Q,Q]=0$ and $Q$ of weight 1, on $\Ms_E$. }\end{ex}

The {\it graded tangent bundle} $\T\Ms$ of an $n$-graded manifold $\Ms$ is by definition the $(n+1)$-graded
manifold associated with the tangent bundle $\sT(\text{\rm dgr}(\Ms))$ of the $n$-vector bundle $\text{\rm
dgr}(\Ms)$ being the degradation of $\Ms$, i.e. $\Ms=\Ms_{\dgr(\Ms)}$. The degree of $\pa_{x^j}$ as a function on
$\dts\Ms$ is $-g(x^j)\in-\N^{n}$. Similarly, to obtain an $(n+1)$-vector bundle associated with $\dts\text{\rm
dgr}(\Ms)$ -- the {graded cotangent bundle} $\T^*\Ms$ -- we have to define the degree of $\pa_{x^j}$ as
$(\ao-g(x^j),1)\in\N^{n+1}$. One can also say that the grading in $\T^*\Ms$ is induced by the Euler vector
fields $\zD_{\T^*\Ms}$ and the phase lifts $\T^*(\zD^k_\Ms)$ which, in the standard adapted local coordinates
$(x^j,p_j)$ have the form $\zD_{\T^*\Ms}=\sum_jp_j\pa_{p_j}$ and
$$\T^*(\zD^k_\Ms)=\sum_j\left(g_k(x^j)x^j\pa_{x^j}+(1-g_k(x^j)p_j\pa_{p_j}\right).$$
The bases of the corresponding projections are $\Ms$ and $\Ms^*_{(k)}$, $k=1,\dots,n$, and this set of
$n$-graded manifolds is closed with respect to the corresponding dualities.

We say that an $r$-form (resp., an $r$-vector field) $\za$ is of weight $i\in\Z^n$, $w(\za)=i$, if
$\Ll_{\zD^k_\Ms}(\za)=i_k\za$, $k=1,\dots,n$. Note that with this convention the weight of $\xd x^j$ is
$w(x^j)=g(x^j)\in\N^n$, but the degree of $\xd x^j$ as a function in $\T\Ms$ is $(g(x^j),1)\in\N^{n+1}$.
Similarly, the weight of $\pa_{x^j}$ is $-w(x^j)=-g(x^j)\in\Z^n$, but the degree of $\pa_{x^j}$ as a function
in $\T^*\Ms$ is $(\ao-g(x^j),1)\in\N^{n+1}$.

\section{Multi-graded symplectic and Poisson manifolds}
\begin{definition}{\rm A {\it $n$-graded symplectic (resp. Poisson) manifold} is an $n$-graded manifold
equipped with a symplectic form of weight $\ao$ (resp., a Poisson tensor of weight $-\ao$).}
\end{definition}
\noindent Recall that a differential $2$-form $\omega$  can be locally written in local coordinates $(x^i)$ as
\begin{equation}\label{e:om}
\om= \frac{1}{2}\sum_{i,j}  dx^i \,\om_{i,j}\,dx^j.
\end{equation}
A $2$-from $\om$ on $\Ms$ is called symplectic, if $d\om = 0$ and $\om$ is
non-degenerate. The latter means that the induced homomorphism of $\As(\Ms)$-modules
\[
\tilde{\om}:\Gamma(\T\Ms)\to \Gamma(\T^*\Ms), \quad X\mapsto i_X\om,
\]
is invertible. For any $(n-1)$-graded manifold $\Ms$, the $n$-graded manifold $\T^*\Ms$ posses a canonical
symplectic form $\zw_\Ms$ of weight $\ao$. Indeed, fixing local coordinates $(x^j)$ in $\Ms$, one can easily
see that the 2-form $\zw_\Ms$ which in the adapted coordinates $(x^j,p_s)$ in $\T^*\Ms$ reads
$$\zw_\Ms=\sum_j\xd p_j\xd x^j$$ is well defined, symplectic and, since $w(\xd p_j)=\ao-w(\xd x^j)$, of weight $\ao$.
Note that for any vector field $X$ on $\Ms$ we can define its phase lift $\T^* X=\dTs X+\zD_{\T^*\Ms}$ exactly
like in the standard case. Here, $\dTs X$ is the cotangent lift of $X$ -- the hamiltonian vector field of the
linear function on $\T^*\Ms$ associated with $X$. The phase lifts are linear, i.e. commute with
$\zD_{\T^*\Ms}$ and satisfy $\Ll_{\T^* X}\zw_\Ms=\zw_\Ms$. We have also full analogs of Propositions \ref{p3}
and \ref{p3a}.

It is well known that any symplectic vector bundle $(E,\zD_E,\zW)$, i.e. a vector
bundle $(E,\zD_E)$ equipped with a symplectic form which is 1-homogeneous with respect
to the Euler vector field, $\Ll_{\zD_E}\zW=\zW$, is canonically isomorphic to the
cotangent bundle over the base of $E$ with the canonical symplectic form:
$(E,\zD_E,\zW)\simeq(\sT^*M,\zD_{\sT^*M},\zw_M)$. A similar fact holds for symplectic
N-manifolds of degree 1 in the terminology of D.~Roytenberg or {\it 1-graded symplectic
manifolds} in our terminology: every 1-graded manifold $\Ms$ equipped with a symplectic
form $\zW$ of weight 1 ($\Ll_{\zD_\Ms}\zW=\zW$) is diffeomorphic to $\T^* M$ equipped
with the canonical symplectic form (cf. \cite[Proposition 3.1]{Roy}). We can generalize
this fact, i.e. we have the following graded version of Theorem \ref{tsym}.
\begin{theo}\label{tgsym} Any $n$-graded symplectic manifold is canonically
isomorphic to the graded cotangent bundle $\T^*\Ms$ of an $(n-1)$-graded manifold $\Ms$, equipped with the
canonical symplectic form $\zw_\Ms$. Moreover, we have canonical symplectomorphisms
\be{sp}(\T^*\Ms,\zw_\Ms)\simeq(\T^*\Ms^*_{(k)},\zw_{\Ms^*_{(k)}}).\ee
\end{theo}
\bepf The proof is completely parallel to that of Theorem
\ref{tsym} and we omit it. \epf

\medskip
Recall (cf. \cite{GM1,GM2}) that a graded Poisson bracket of degree  $i$  on  a  $\Z^n$-graded associative
commutative algebra $\A=\oplus_{k\in\Z^n}\A^k$ is a graded bilinear map $ \{\cdot,\cdot\}:\As\ti\As\ra\As$ of
degree $i\in\Z^n$ such that
\begin{enumerate}
\item $\{ a,b\}=-(-1)^{(\vert a\vert+\vert i\vert)(\vert b\vert+\vert i\vert)}\{ b,a\}$ (graded
anticommutativity), \item $\{ a,bc\}=\{ a,b\} c+(-1)^{(\vert a\vert+\vert i\vert)\vert b\vert}b\{  a,c\}$
(graded Leibniz rule), \item $\{\{ a,b\},c\}=\{ a,\{ b,c\}\}-(-1)^{(\vert a\vert+\vert i\vert)(\vert
b\vert+\vert i\vert)}\{ b,\{ a,c\}\}$ (graded Jacobi identity), where $\vert a\vert$ denotes the total degree
of $a$, etc.
\end{enumerate}
A homogeneous element $q$ of degree $k$ with the parity opposite to the parity of $i$ we call {\it
homological} if $\{ q,q\}=0$. It induces a cohomology operator $\xd_q=\{ q,\cdot\}$ of odd total degree $\vert
k+i\vert$ on $\A$.

The bracket $\{\cdot,\cdot\}=\{\cdot,\cdot\}_\Ms$ on $\As(\T^*\Ms)$ associated with the canonical symplectic
form $\zw_\Ms$ and represented locally by the Poisson tensor $\zL_\Ms=\sum_j\pa_{p_j}\pa_{x^j}$ is a graded
Poisson bracket of degree $-\ao$. Since the algebra $\As(\T^*\Ms)$ is non-negatively graded, negative degrees
mean simply 0. Recall also (cf. \cite{Roy0,Roy})that with any linear Poisson structure $\zL$ on a vector
bundle $E$, thus with any de Rham differential $D_\zL$ on the Grassmann algebra $A(E^*)$ of multi-sections of the
dual bundle, one can associate a function $H_\zL=\zi_{D_\zL}$ on $\T^*\Ms_E$, so that $\{ H_\zL,\cdot\}$ is the cotangent lift $\dTs D_\zL$.
Here, we clearly identify $\As(\Ms_E)=A(E^*)$ with basic functions on $\T^*\Ms_E$.
In local homogeneous coordinates $(x,y)$ on $E$ and $(x,\zvy,p,\zx)$ on $\T^*\Ms_E$, every $\zL$ of weight -1 is of the
form
$$\zL=\sum_{a,r}\zr^r_a(x)\pa_{y^r}\we\pa_{x^a}+\frac{1}{2}\sum_{r,s}C^u_{r,s}(x)y^u\pa_{y^r}\we\pa_{y^s},$$
so
$$H_\zL=\sum_{a,r}\zr^r_a(x){\zx^r}{p_a}-\frac{1}{2}\sum_{r,s}C^u_{r,s}(x)\zvy_u\zx^r\zx^s\,.$$
In other words the function $H_\zL$ is a homological Hamiltonian whose hamiltonian vector field $Q_\zL=\{
H_\zL,\cdot\}$ is the cotangent lift of the de Rham derivative $D_\zL$.
As $\dts\Ms_E\simeq\dts\Ms_{E^*}$,
this Hamiltonian can be also viewed
as the function $\zi_\zL$ on $\dts\Ms_E$ associated with the Poisson tensor $\zL\in A(E)$. The Hamiltonian
$H_\zL$ has not only total degree 3 but it is homogenous of bi-degree $(1,2)$. In some terminology one says
also that $\zL$ determines a {\it Lie algebroid structure} on $E$ and that $(\Ms,Q_\zL)$ is a {\it
$Q$-manifold} or {\it Lie antialgebroid} in the language of \cite{AKSZ,Va,Vor1,Vor2}. To find a generalization
for Poisson $n$-bundles $(F,\zL)$ let us recall that in this case the Poisson tensor $\zL$ determines
(concordant) Lie algebroid structures on all vector bundles $F^*_{(k)}\ra F_{[k]}$, i.e. (by definition
{\it concordant}) homological vector fields $q_{[k]}$ of degree $\zd^n$ on $\Ms_{F^*_{(k)}}$. We will say that
a vector field on an $n$-graded manifold is {\it unital} if its homogeneous parts have weights
$\zd^1=(1,0,\dots,0)$, $\zd^2=(0,1,0,\dots,0)$, etc. Observe that any unital vector fields $q_{[k]}$ on the
side bundle $\Ns_{[k]}$ of an $(n+1)$-graded symplectic manifold $\Ns$ defines the {\it induced} vector field
$(q_{[k]})_{[s]}$ on $\Ns_{[s]}$, $s=1,\dots,n+1$, defined as the restriction of the cotangent lift $\dTs
q_{[k]}$ to $\T^*\Ns_{[k]}\simeq\Ns$. If we consider the $(n+1)$-graded symplectic manifold $\Ns=\T^*\Ms_F$,
then $\Ms_{F^*_{(k)}}=\Ns_{[k]}$. We have the following.
\begin{prop}\label{conc} Homological vector fields $q_{[k]}$ of degree $\zd^n$ on $\Ms_{F^*_{(k)}}=\Ns_{[k]}$,
$k=1,\dots,n$, are concordant if and only if their cotangent lifts coincide (up to the identification
$\T^*\Ms_{F^*_{(k)}}\simeq\T^*\Ms_{F^*_{(s)}}$), i.e., if and only if $(q_{[k]})_{[s]}=q_{[s]}$ for all
$k,s=1,\dots,n$.
\end{prop}
\bepf
The cotangent lifts, uniquely determined by their restrictions to $\Ms_{F^*_{(k)}}=\Ns_{[k]}$, can be easily
seen as represented by the Hamiltonian vector field with the Hamiltonian $H_\zL$ associated with $\zL$. \epf

\section{Higher Courant structures, higher Lie algebroids, and Drinfeld $n$-tuples}
A {\it Lie bialgebroid}, as introduced in \cite{MX}, is a pair of Lie algebroid structures $\zL,\zL^*$ on a
vector bundle $E$ and its dual $E^*$ that satisfies certain compatibility condition. This compatibility
condition has been recognized in \cite{Roy0} as the commutation of the corresponding homological Hamiltonians
$\{ H_\zL,H_{\zL'}\}=0$ on $\T^*\Ms_E\simeq\T^*\Ms_{E^*}$. This means exactly that the Hamiltonian
$H=H_\zL+H_{\zL'}$ of total degree 3 is homological and concentrated in bi-degrees $(1,2)+(2,1)$, i.e. the
corresponding hamiltonian vector field $Q$ is concentrated in weights $(1,0)+(0,1)$. The total weight of $Q$
is 1. Note that there are homological vector fields of total weight 1 and bi-degrees $(-1,2)$ or $(2,-1)$. They lead
to the concept of {\it quasi Lie bialgebroids}. The {\it derived bracket} (in the terminology of
Y.~Kosmann-Schwarzbach)
$$\{ X,Y\}_Q=\{\{ X,H\},Y\}=-(-1)^{\vert x\vert+n}\{ Q(X),Y\},$$
is closed on functions representing sections of $E\ti_ME^*$ and gives a standard model of a {\it Courant
bracket} \cite{Co,Do} in its non-symmetric version, or a {\it Courant algebroid} \cite{Co,LWX,Roy0}. Note
however that the concept of Courant algebroid is more general and based on graded symplectic manifolds of
degree 2 which are not bi-graded in general. The whole structure, i.e. $\T^*\Ms_E$ with the canonical
symplectic Poisson bracket and the homological hamiltonian and unital vector field $Q$ -- the {\it Drinfeld
double} of the original Lie bialgebroid -- is a natural generalization of the Drinfeld double Lie algebra
\cite{Dri}.

A natural generalization of the above concepts is as follows.
\begin{definition}{\rm An {\it $n$-Courant structure} is an $n$-graded symplectic manifold $(\Ns,\zW)$
equipped with a homological Hamiltonian of total degree $(n+1)$. A {\it Drinfeld
$n$-tuple} is an $n$-Courant structure whose homological hamiltonian vector field is
unital.}
\end{definition}\noindent
We know already that $(\Ns,\zW)=(\T^*\Ms,\zw_\Ms)$ for an $(n-1)$-graded manifold
$\Ms$. Let $\{\cdot,\cdot\}$ be the symplectic Poisson bracket on
$(\Ns,\zW)=(\T^*\Ms,\zw_\Ms)$. This Poisson bracket is a graded Poisson bracket on
$\As(\Ns)$ of degree $-\ao$, i.e. of total degree $-n$. It is a general algebraic fact
that the {\it derived bracket} $\{ X,Y\}_H=\{\{ X,H\},Y\}$ of any homological
Hamiltonian $H\in\As^{n+1}(\Ns)$ of the total degree $(n+1)$ (or of any homological
hamiltonian and vector field $Q$ of total weight 1) is then a {\it Leibniz bracket} of
total degree $(1-n)$. In classical terms, this bracket can be interpreted as a {\it
$n$-Courant algebroid}, i.e. as a bracket on $\As^{(n-1)}(\Ns)$ -- the module of
sections of the vector bundle $\Cs=V^{(n-1)}(\Ns)$ over $M$. Note that the symplectic
Poisson bracket also gives rise to certain operations on subbundles of
$\Cs=V^{(n-1)}(\Ns)$, namely $\la\cdot,\cdot\ran_{j,k}:V^{j}(\Ns)\otimes_M V^k(\Ns)\ra
V^{j+k-n}(\Ns)$, or, globally, to a graded operation
$$\la\cdot,\cdot\ran:\Cs(\Ns)\otimes_M\Cs(\Ns)\ra \Cs(\Ns)$$
of degree $-\ao$ on the graded vector bundle
$$\Cs(\Ns)={\bigoplus_{\vert i\vert<n}}{}_MV^i(\Ns).$$
If we do not insist on working with multi-graded symplectic manifolds, which means -- cotangent bundles -- and
we admit symplectic graded N-manifolds of degree $n$ in the terminology of D.~Roytenberg \cite{Roy}, then,
analogously, we get a notion of a {\it Courant algebroid of degree $n$}. In this context, however, the concept
of Drinfeld $n$-tuple makes no sense.

Recently, the double Lie algebroids, as introduced by K.~C.~Mackenzie \cite{Mac1a} -- \cite{Mac1e}, have been
recognized by T.~Voronov \cite{Vor2} as double $Q$-manifolds and generalized to $n$-fold $Q$-manifolds, i.e.
$n$-graded manifolds $\Ms_F$ (associated with an $n$-vector bundle $F$) and endowed with a homological unital
vector field $Q$. This means that $Q=Q_1+\cdots+Q_n$, where $Q_1,\dots,Q_n$ are commuting homological vector
fields of $n$-degrees, respectively, $\zd^1,\dots,\zd^n$.

More precisely, an {\it $n$-fold Lie algebroid} is an $n$-vector bundle $F$ equipped with Lie algebroids
structures on the vector bundles $h^k_0:F\ra F_{[k]}$, $k=1,\dots,n$, and satisfying certain compatibility
conditions. In particular, all morphism in the characteristic diagram should be Lie algebroid morphisms. An
elegant way to describe these conditions is to pass to the corresponding $n$-graded manifold $\Ms=\Ms_F$.
Then, we can interpret these Lie algebroid structures as homological vector fields $Q_k$ of weight 1 on the
corresponding 1-graded manifolds and finally, due to the commutativity of Lie algebroid morphisms, as a
homological and unital vector field on the $n$-graded total space $\Ms$. The compatibility conditions reduce
now to the fact that the vector fields $Q_k$ commute. Equivalently, the total vector field $Q=Q_1+\cdots+Q_n$
is unital and homological, so we end up with the following (see \cite{Vor2}).
\begin{definition}{\rm An {\it $n$-fold Lie algebroid} is an $n$-graded manifold with a homological and unital
vector field.}
\end{definition}\noindent
Observe that $n$-fold Lie algebroid is a particular case of a Drinfeld $(n+1)$-tuple.
\begin{prop} There is a one-to-one correspondence between $n$-fold Lie algebroids $(\Ms,Q)$
and $(n+1)$-Drinfeld tuples $(\dts\Ms,\wt{Q})$ such that $\wt{Q}$ has trivial summand of weight $\zd^{n+1}$.
\end{prop}
\bepf Let us put $\wt{Q}=\dTs Q$. Then, as easily seen, $\wt{Q}$ is a Hamiltonian and homological vector field on $\dts\Ms$ of
weight $(0,w(Q))$. Conversely, if $\wt{Q}$ is a Hamiltonian and unital vector field on $\dts\Ms$ with the
trivial summand of weight $\zd^{n+1}$, then $\wt{Q}=\dTs Q$ for some unital vector field $Q$ on $\Ms$.
Moreover, since $\dTs[Q,Q]=[\wt{Q},\wt{Q}]=0$, then $Q$ is homological. \epf

\medskip
\noindent
%%%%%%%%%%%%%%%%%% MR %%%%%%%%%%%%%%%%%%%%%%%%%%
Let $F$, $\Ms_F$, $Q$ be as above and let $i,j\in\{0,1\}^n$, $i<j$. We claim that the vector field $Q$ induces
a $(|j|-|i|)$-fold Lie algebroid structure on the $(|j|-|i|)$-vector bundle whose total space is $F_j$ and the final
base is $F_i$.
%Let us denote this
%bundle by $_iF_j$ and the corresponding $(|j|-|i|)$-graded
%manifold by $_i\Ms_j$.
Because $Q$ is tangent to $\Ms_j$ it is enough to verify the claim for $\Ms = \Ms_F$, i.e. for $j=\ao$.
Moreover, an inductive reasoning shows that we may also assume that $|i|=1$, since an $(n-|i|)$-vector bundle
can be reached in $|i|$ steps in which we simply forget about one Euler vector field. Note that the
$(n-1)$-graded manifold $_{\{k\}}\Ms$ associated with the $(n-1)$-vector bundle ${_{\{k\}}F}$ defined by Euler
vector fields $\{\zD^s\}$, $s\neq k$, can be obtained from $\Ms_F$ by applying the parity changing functor to
the super vector bundle $\Ms\to\Ms_{[k]}$, i.e. $_{\{k\}}\Ms = \Pi_{\Ms_{[k]}}\Ms$. The vector field
$\tilde{Q} := Q - Q_k$ is a linear vector field with respect to the vector bundle $\Ms\to \Ms_{[k]}$. The
following general fact implies that we can pass from $\tilde{Q}$ to a vector field $\bar{Q}$ on $_{\{k\}}\Ms$.
\begin{lem} Let $\Es\to\Ms$ be a super vector bundle and $Vect_0(\Es)$
denotes the Lie algebra of linear vector fields on $\Es$. There exist a canonical $\As(\Ms)$-linear
isomorphism of Lie algebras
\be{isoL} \phi: Vect_0(\Es)\to Vect_0(\Pi_{\Ms}\Es). \ee
\end{lem}
\bepf Let $\{\eta_i\}$ and $\{\mu_i\}$ be the corresponding local linear coordinates on $\Es$ and
$\Pi_{\Ms}\Es$, respectively, and let $\tilde{i} := \tilde{\eta}_i$, so $\tilde{\mu}_i = \tilde{i}+1$. The Lie
algebra $Vect_0(\Es)$ is locally spanned by the vector fields $\eta_i\pa_{\eta_j}$. The following formula
\[
\phi(\eta_i\pa_{\eta_j}) = (-1)^{\tilde{i} + \tilde{j}} \mu_i\pa_{\mu_j}
\]
does not depend on the choice of local coordinates. In fact, if $\eta_i'=\sum_j \eta_j T_{ji}(x)$, $T_{ji}(x)
\in \As(\Ms)$, describes transformations of fiber coordinates of $\Es$, then also $\mu_i'=\sum_j \mu_j
T_{ji}(x)$ and
\beas
\eta_i'\pa_{\eta_j'} &=& \sum_k \eta_i'\frac{\pa\eta_k}{\pa\eta_j'}\pa_{\eta_k}  =
\sum_{l,k} \eta_lT_{li}(x) T^{jk}(x) \pa_{\eta_k}\\
&=& \sum_{l,k} (-1)^{\tilde{l}(\tilde{i} + \tilde{l} + \tilde{k} +\tilde{j})}
T_{il}(x)T^{jk}\eta_l\pa_{\eta_k},
\eeas
because the parity of $T_{ij}(x)$ is $\tilde{i}+\tilde{j}$. Hence
\[
\phi(\eta_i'\pa_{\eta_j'})  = \sum_{l,k} (-1)^{\tilde{i} + \tilde{j}}  \mu_lT_{il}(x)T^{jk}(x)\pa_{\mu_k} =
(-1)^{\tilde{i}+ \tilde{j}} \mu_i'\pa_{\mu_j'}.
\]
It is also easy to calculate that $\phi$ preserves the Lie bracket of vector fields. \epf

\medskip
\noindent It follows from above lemma that $\bar{Q}$ is a homological vector field on $_{\{k\}}\Ms$ and so
induces a $(n-1)$-fold Lie algebroid structure on $_{\{k\}}F$.
%%%%%%%%%%%%%%%%%%%
Thus we get the following.
\begin{prop} An $n$-fold Lie algebroid structure on an $n$-vector bundle $F$ induces canonically, for $i,j\in\{ 0,1\}^n$,
$i<j$, an $(|j|-|i|)$-fold Lie algebroid structure on the $(|j|-|i|)$-vector bundle whose total space is $F_j$
and the final base is $F_i$.
\end{prop}\noindent
%%%%%%%%%%%%%%%%%%MR%%%%%%%%%%%%%%%%%%%%%%%%%%%%
A natural way of constructing $n$-fold Lie algebroids can be based on the following trivial observation.
\begin{theo}\label{txx}
If $(\Ms,Q)$ is an $n$-fold Lie algebroid, then $(\T\Ms,\dt Q+\xd)$, where $\xd$ is the de Rham differential
on $\T\Ms$, is an $(n+1)$-fold Lie algebroid.
\end{theo}
\bepf
Since the tangent lift of vector fields respects the Schouten brackets \cite{GU}, $\dt Q$ is a homological
vector field with components of weights $\zd^1,\dots,\zd^n$. Moreover, any tangent lift, which locally reads
$$\sum_a\left(f^a(x)\pa_{x^a}+\sum_b\frac{\pa f^a}{\pa x^b}\dot{x}^b\pa_{\dot{x}^a}\right),$$ commutes with the de Rham
vector field $\xd=\sum_a\dot{x}^a\pa_{x^a}$ of weight $\zd^{n+1}$, so $Q+\dt$ is homological and unital.
 \epf

\medskip\noindent
The question now is: what is the higher analogue of a Lie bialgebroid? Our answer is obvious: it corresponds
to a Drinfeld $n$-tuple. Recall that a  Drinfeld $n$-tuple is an $n$-graded symplectic manifold with a
homological Hamiltonian and unital vector field. Let us take an $n$-graded symplectic manifold $(\Ns,\zW)$. We
know already that $(\Ns,\zW)=(\T^*\Ms,\zw_\Ms)$ for an $(n-1)$-graded manifold $\Ms$. If $Q=Q_1+\cdots+Q_n$ is
the decomposition of the homological vector field of a Drinfeld $n$-tuple on $(\Ns,\zW)=(\T^*\Ms,\zw_\Ms)$
into homogeneous parts of weights $\zd^1,\dots,\zd^n$. It is easy to see that $[Q,Q]=0$ implies $[Q_r,Q_s]=0$
for all $r,s=1,\dots,n$, i.e. all vector fields $Q_k$ are homological and pair-wise commuting. Moreover, for
any $s\ne r$, the vector field $Q_r$ is tangent to $\Ns_{[s]}$  and projectable to the vector field
$q^r_{[s]}$ with respect to the canonical projection $\Ns\ra\Ns_{[s]}$.

Therefore, $\Ns_{[s]}$ posses canonically $n$ commuting homological vector fields $q^r_{[s]}$, $r,k=1,\dots,n$
-- the restrictions of $Q_1,\dots,Q_n$. Here, for technical convenience, we list $n$ vector fields for each
base, but clearly $q^r_{[r]}=0$, i.e. $\Ns_{[s]}$ is canonically an $(n-1)$-fold Lie algebroid. But the
collection of all $\Ns_{[s]}$ is a collection of $(n-1)$-graded manifolds
$\Ms,\Ms^*_{(1)},\dots\Ms^*_{(n-1)}$, closed with respect to duality. As a matter of fact, according to the
graded analog of Proposition \ref{p3}, the vector field $Q_r$ is the cotangent lift of $q^r_{[s]}$ for $s\ne
r$. The compatibility condition for all the $(n-1)$-fold Lie algebroid structures on
$\Ms,\Ms^*_{(1)},\dots\Ms^*_{(n-1)}$ is expressed by saying that they come from projections of certain
homological hamiltonian and unital vector fields on $\T^*\Ms$, i.e. from Drinfeld $n$-tuples. We can say
that $n$-tuple Lie algebroid corresponds to a Drinfeld $n$-tuple, exactly like a Lie bialgebra (or a Lie
bialgebroid) corresponds to a Drinfeld double Lie algebra (or Lie algebroid). In particular, these structures
are compatible with the base structure $\Bs\Ns$, i.e. the projections of $q^r_{[k]}$ and $q^r_{[s]}$ on
$\Ns_{[k,s]}$ coincide, where $[k,s]\in\{ 0,1\}^{n}$ has zeros exactly at positions $k,s=1,\dots,n$. The pair
of homological vector fields: $q^k_{[s]}$ on $\Ns_{[s]}$, and $q^s_{[k]}$ on $\Ns_{[k]}$, where $E=\Ns_{[s]}$
and $E^*=\Ns_{[k]}$ are regarded as dual vector bundles over $N_{[k,s]}$ -- the degraded manifold
$\Ns_{[k,s]}$, forms a Lie bialgebroid. Indeed, $Q_k$ and $Q_s$ are the cotangent lifts of $q^k_{[s]}$ and
$q^s_{[k]}$ to $\Ns\simeq\T^*\Ns_{[s]}\simeq\T^*\Ns_{[k]}$, so in $\dts E\simeq\dts E^*$ they are represented
by commuting and homological Hamiltonians $H_k$ and $H_s$ of degrees $(2,1)$ and $(1,2)$. According to the
result of D.~Roytenberg \cite{Roy0}, this means exactly that we deal with a Lie bialgebroid. It is not true,
however, that in general all these Lie bialgebroid structures produce a Drinfeld $n$-tuple even under a
natural condition saying that the vector fields $q^r_{[k]}$ and $q^r_{[l]}$ induced from algebroid structures
on $N_{[k]}\to N_{[r,k]}$ and $N_{[l]}\to N_{[r,l]}$ coincide on $\Ns_{[k,l]}$, as shown in the following example.

\begin{ex} {\rm Consider a trivial double vector bundle $M = \R\times\R^2\times \R\times \{*\}$ over a point
$\{*\}$ with $1$-dimensional core
and trivial side bundles of rank $1$ and $2$. Then $\Ns :=\T^*M$ carries a $3$-vector bundle structure. Let us
denote by $v_{011}$,  $v_{101}$ and $v_{001}^{(1)}, v_{001}^{(2)}$ the fiber coordinates on the core and the
side vector bundles $N_{[1,3]}$ and $N_{[1,2]}$, respectively. Let $\{e_{001}^{(1)}, e_{001}^{(2)}, e_{011}\}$
be the corresponding dual basis of sections of the bundle $N_{[1]} \to N_{[1,3]}$. We endow this bundle with a
structure of Lie algebroid by setting the anchor to zero and $\As(N_{[1,3]})$-linear Lie bracket as follows:
\[
[e_{001}^{(1)}, e_{001}^{(2)}]:= v_{010}\cdot e_{011},
\]
$[e_{001}^{(k)}, e_{011}] :=0$ for $k=1,2$. Obviously, this is a nilpotent Lie bracket. The induced
homological vector field $q^3_{[1]}$ on $\Ns_{[1]}$ in the corresponding graded local coordinates
$\{\theta_i^a\}$, $i\in \{0,1\}^3$, has the form
\[
q^3_{[1]} = \theta_{010}\theta_{001}^{(1)}\theta_{001}^{(2)}\pa_{\theta_{011}}.
\]
Let us assume that the other $5$  vector bundles $N_{[k]}\to N_{[k,l]}$ carries the zero Lie algebroid
structure. Then of course the $3$ pairs of Lie algebroids $(N_{[k]}, N_{[l]})$ over $N_{[k,l]}$ constitute a
Lie bialgebroid. Note that the restriction of $q^3_{[1]}$ to $N_{[1,2]}$ is zero, so the vector fields
$q^r_{[k]}$ coincide on intersections, i.e. $q^r_{[k]}|N_{[k,l]} = q^r_{[l]}|N_{[k,l]} = 0$ for distinct
$r,k,l$. However they do not come from Drinfeld $3$-tuple on $\Ns$, because the cotangent
lifts $\dt^*q^3_{[1]} \neq 0$ and $\dt^*q^3_{[2]} =0$ (as $q^3_{[2]}
=0$)
can not be equal to $Q_3$ at the same time.} 

%because the hamiltonian
%$h^3_{[1]}=\theta_{010}\theta_{001}^{(1)}\theta_{001}^{(2)}\theta_{100}$ associated with $q^3_{[1]}$
% is not zero, while $h^3_{[2]}=0$. }
\end{ex}

\medskip\noindent
The following theorem gives sufficient conditions.
%%%%%%%%%%%%%%%%%%%%%%%%%%%%MR%%%%%%%%%%%
\begin{theo}\label{t11} Let $\Ms$ be an $(n-1)$-graded manifold.
A Drinfeld $n$-tuple on $\Ns=\T^*\Ms$ is equivalent to a collection of Lie bialgebroid structures on all the
pairs of dual vector bundles $N_{[k]}$ and $N_{[s]}$ over the common base $N_{[k,s]}$, $k\ne s$, related to
homological unital vector fields $q^s_{[k]}$ and $q^k_{[s]}$, respectively, which satisfy the compatibility
condition
\be{ccc} (q^r_{[k]})_{[s]}=q^r_{[s]}\quad\text{for}\quad r\ne k
\ee
with the convention $q^s_{[s]}=0$.
\end{theo}
\bepf The vector field $q^r_{[k]}$ is of degree 1 with respect to the Euler vector field $\zD^r$, so of degree
0 with respect to $\zD^s$, $s\ne r$. A similar statement is true for the vector field $Q^r_{[k]}$ -- the
(unique) linear Hamiltonian extension of $q^r_{[k]}$ to $\Ns$, i.e. $Q^r_{[k]}=\dTs(q^r_{[k]})$. According to
the graded version of Proposition \ref{p3}, the vector field $Q^r_{[k]}$ coincides with $Q^r_{[s]}$ for
$k,s\ne r$, so that $Q_r$ defined as $Q^r_{[s]}$ for some (thus all) $s\ne r$ is of degree $\zd^r$. Since
$Q^r_{[s]}$ and $Q^s_{[r]}$, $r\ne s$, commute because they correspond to the given bialgebroid structure on
$(N_{[s]}, N_{[r]})$, the vector fields $Q_r$ pairwise commute and $Q=Q_1+\cdots+Q_n$ gives rise to a Drinfeld
$n$-tuple on $\Ns=\T^*\Ms$ which induces prescribed bialgebroid structures. \epf

\medskip\noindent
Since, for a fixed $r$, the vector fields $q^r_{[k]}$, $k\ne r$, are concordant, we can also characterize a
Drinfeld $n$-tuple in terms of Poisson structures.
\begin{theo} A Drinfeld $n$-tuple is a collection of Poisson $(n-1)$-vector bundles in duality: $(F,\zL)$ and
$(F^*_{(k)},\zL_k)$, $k=1,\dots,n-1$, which are compatible in the sense that the corresponding Hamiltonians
$H_\zL$ and $H_{\zL_k}$, $k=1,\dots,n-1$, interpreted as functions on $\T^*\Ms_F\simeq\T^*\Ms_{F^*_{(k)}}$,
commute with respect to the symplectic Poisson bracket.
\end{theo}

\medskip
Let us end up with some words about reduction. Since we deal, in fact, with a homological Hamiltonian system
$(\Ns,\zW,H)$ on a symplectic super-manifold $(\Ns,\zW)$, the reduction should be understood as the
Hamiltonian reduction with respect to a coisotropic and $n$-graded submanifold $\Ns_0$. If we assume that the
Hamiltonian $H$ is constant on leaves of the characteristic foliation $\cF$ of $\Ns_0$ and the quotient
$\Ns'=\Ns_0/\cF$ is a well-defined multi-graded manifold, then the restriction of $\zW$ to $\Ns_0$ projects to
a symplectic form $\zW'$ on $\Ns'$, the homological Hamiltonian $H$ projects to a homological Hamiltonian $H'$
on $\Ns'$ and we end up with a new homological Hamiltonian system $(\Ns',\zW',H')$ -- new $n$-Courant or new
Drinfeld $n$-tuple structure. Of course, this picture covers the reduction associated with the moment map of a
Hamiltonian group action: this is only the choice of the coisotropic submanifold which is determined by the
moment map $\zm$ -- the inverse-image $\zm^{-1}(\{ 0\})$. Note only that this group action should respect the
graded structure, i.e. it should commute with the Euler vector fields.
\begin{ex}{\rm Consider the canonical symplectic triple vector bundle
$$\dts\sT\sT M\simeq\dts\dts\sT M\simeq\dts\sT\dts M\simeq\dts\dts\dts M$$
with the characteristic diagram

\[
\xymatrix{
               &  \dts\sT M \ar[dd]\ar[dl]   &                     &  \dts\sT\sT M \ar[ll] \ar[dd] \ar[dl]
  \\
\dts M \ar[dd] &            & \sT\dts M \ar[ll]\ar[dd] &
  \\
               &  \sT M \ar[ld] &                 &  \sT\sT M \ar[ll] \ar[ld]
  \\
M        &                   &      \sT M\ar[ll] & }
\]

\medskip\noindent
The tangent bundle $\zt_M:\sT M\ra M$ is canonically a Lie algebroid. The corresponding homological vector
fields on $\T M$ is the de Rham vector field $D_M$ which in local coordinates $(x,\dot{x})$ (we do not
distinguish coordinates in $\sT M$ and $\T M$, etc., for simplicity) has the form
$D_M=\sum_a\dot{x}^a\pa_{x^a}$, so the corresponding hamiltonian of degree $(2,1)$ on $\T^*\T M$ is
$H_{(2,1)}=\sum_a\dot{x}^ap_a$. It is well known that Lie algebroid structures
on $\dts M$ such that, together with $D_M$, give a Lie
bialgebroid come from Poisson structures $\zL=\frac{1}{2}\sum_a\zL_{ab}(x)\pa_{x^a}\we\pa_{x^b}$ on $M$. The
corresponding Hamiltonian on $\T^*\T^* M\simeq \T^*\T M$ of degree $(1,2)$ is associated with the tangent lift
$\dt \zL$ of $\zL$ by $H_{(1,2)}=\zi_{\dt\zL}$, i.e.
$$H_{(1,2)}=\sum_{a,b}\zL_{ab}(x){p_a}\dot{p}_b+\frac{1}{2}\sum_{a,b,c}\frac{\pa\zL_{ab}}{\pa x^c}(x)
\dot{x}^c\dot{p}_b\dot{p}_a\,.$$ Take now a Lie group action $G\ti M\ra M$ which is free and proper, so that
the space of orbits $M/G$ is a manifold, and which preserves $\zL$, so that $\zL$ projects onto a Poisson
structure $\zL'$ on $M/G$. Let $(y^s)$ be a basis in the Lie algebra $\cal G$ of $G$ and let
$Y^s=\sum_af^s_a(x)\pa_{x^a}$ be the corresponding fundamental vector fields of this action. Preserving $\zL$
by the action means that the Schouten brackets $[Y^s,\zL]$ vanish. By means of the tangent functor, the action
of $G$ can be extended to an action of the group $\sT G$ on $\sT M$. The Lie algebra of $\sT G$ is $\sT \cal
G$ and the fundamental vector fields of this action are the tangent and vertical lifts, $\dt Y^s$ and $\vt
Y^s$, of the fundamental vector fields of the action of $G$ (cf. \cite{GU}). Since $\dt$ respects the Schouten
bracket and $[\vt Y,\dt\zL]=\vt[Y,\zL]$ (see \cite{GU}), the extended action preserves $\dt\zL$. Moreover, it
is easy to see that $\sT M/\sT G\simeq\sT(M/G)$ and that the canonical projection $\sT M\ra\sT M/\sT
G\simeq\sT(M/G)$ is a Poisson map of $\dt\zL$ onto $\dt(\zL')$.

Consider now the phase prolongation of the $\sT G$ action, $\sT G\ti\dts\sT M\ra\dts\sT M$. It is a
Hamiltonian action with a canonical equivariant moment map $\zm:\dts\sT M\ra(\sT\cal G)^*$. The Hamiltonians
associated with $\dt Y^s$ and $\vt Y^s$ are, respectively, $\zi_{\dt Y^s}$ and $\zi_{\vt Y^s}$ so that the
submanifold $N_0$ consisting of common zeros of all functions $\zi_{\dt Y^s}$ and $\zi_{\vt Y^s}$, is a
coisotropic submanifold of $\dts\sT M$. Since these functions are linear, it is a vector subbundle of $\dts\sT
M\ra\sT M$. But $\dts\sT M$ is canonically a symplectic double vector bundle with the other projection onto
$\dts M$ and $N_0$ is a vector subbundle also with respect to the other bundle structure. This is because the
tangent lifts of vector fields are linear and the vertical lifts are homogeneous of degree -1, so
$\dTs\zD_{\sT M}(\zi_{\dt Y^s})=0$ and $\dTs\zD_{\sT M}(\zi_{\vt Y^s})=-\zi_{\vt Y^s}$, thus the corresponding
homotheties do not leave $N_0$. The manifold $N_0$ has therefore its graded counterpart $\Ns_0$ being a
bigraded coisotropic submanifold of $\T^*\T M$. The characteristic distribution on $\Ns_0$ is spanned by the
(super) vector fields $\dTs\dt Y^s$ and $\dTs\vt Y^s$. They preserve the Hamiltonian $H=H_{(2,1)}+H_{(1,2)}$
associated with the Lie bialgebroid structure $(D_M,D_\zL)$ and the Hamiltonian reduction leads to the
bi-graded symplectic manifold $\T^*\T (M/G)$ with the reduced homological hamiltonian $H'$ associated with the
Lie bialgebroid structure $(D_{M/G},D_{\zL'})$.

We can go further to the iterated tangent bundle $\sT\sT M$ which is a double vector bundle with respect to
projections $\zt_{\sT M}$ and $\sT\zt_M$ onto $TM$. It is also canonically a double Lie algebroid
corresponding to the homological vector field $q=D_{\sT M}+\dt(D_M)$ on $\T\T M$ which in local coordinates
$(x,\dot{x},\bar{x},\ddot{x})$ takes the form
$$q=\sum_a\left(\bar{x}^a\pa_{x^a}+\ddot{x}^a\pa_{\dot{x}^a}\right)+\sum_a\left(\dot{x}^a\pa_{x^a}+
\ddot{x}^a\pa_{\bar{x}^a}\right).$$ It corresponds to the linear function $\zi_q$ on $\T^*\T\T M$ which in
the adapted local coordinates $(x,\dot{x},\bar{x},\ddot{x},p,\dot{p},\bar{p},\ddot{p})$ of degrees,
respectively,
$$(0,0,0),(1,0,0),(0,1,0),(1,1,0),(1,1,1),(0,1,1),(1,0,1),(0,0,1),$$
reads
$$\zi_q=H_{(1,2,1)}+H_{(2,1,1)}=\sum_a\left(\bar{x}^ap_a+\ddot{x}^a\dot{p}_a\right)
+\sum_a\left(\dot{x}^ap_a+\ddot{x}^a{\bar{p}_a}\right).$$ A homological Hamiltonian $H_{(1,1,2)}$ of degree
$(1,1,2)$ can be obtained from the iterated tangent lift $\dt\dt\zL$  which is linear with respect to both
vector bundle structures:
\beas H_{(2,1,1)}&=&\sum_{a,b}\zL_{ab}(x)\left(\ddot{p}_bp_a+\dot{p}_b\bar{p}_a\right)+
\sum_{a,b,c}\left(\frac{\pa\zL_{ab}}{\pa x^c}(x)\dot{x}^c\dot{p}_a\ddot{p}_b+\frac{\pa\zL_{ab}}{\pa
x^c}(x)\bar{x}^c\bar{p}_a\ddot{p}_b\right)+\\
&&\sum_{a,b,c}\frac{1}{2}\left(\frac{\pa\zL_{ab}}{\pa x^c}(x)\ddot{x}^c+\sum_d\frac{\pa^2\zL_{ab}}{\pa x^c\pa
x^d}(x)\bar{x}^d\dot{x}^c\right)\ddot{p}_b\ddot{p}_a\,.\eeas It clearly commutes with $\zi_q$, so
$H=H_{(2,1,1)}+H_{(1,2,1)}+H_{(1,1,2)}$ represents a Drinfeld triple. The corresponding double Lie algebroid
structures on the side bundles are: $(\T\T M,D_{\sT M}+\dt D_M)$, $(\T\T^* M,\dt D_\zL+D_{\dts M})$, and
$(\T^*\T M,\dTs D_M+D_{\dt D_\zL})$.

Extending the tangent lift action $G\ti\sT M\ra\sT M$ to the iterated tangent lift
action $\sT\sT G\ti\sT\sT M\ra\sT\sT M$ and taking its phase prolongation $\sT\sT
G\ti\dts\sT\sT M\ra\dts\sT\sT M$ (which is canonically Hamiltonian), we get a momentum
map $\zm_1:\dts\sT\sT M\ra(\sT\sT\cal G)^*$ and the corresponding coisotropic
submanifold $N_1=\zm_1^{-1}(\{ 0\})$. This submanifold is 3-homogeneous with respect to
the triple vector bundle structure on $\dts\sT\sT M\simeq\dts\dts\sT M\simeq\dts\sT\dts
M$, so it has its  graded counterpart $\Ns_1$ in the 3-graded symplectic manifold
$\T^*\T\T M\simeq\T^*\T^*\T M\simeq\T^*\T\T^* M$. Like above, the homological
Hamiltonian $H$ is constant on the characteristic distribution of the coisotropic
manifold $\Ns_1$ and we get a reduction to the 3-graded symplectic manifold $\T^*\T\T
(M/G)$ with the reduced homological hamiltonian $H'$ associated with the three
homological vector fields $\dt\dt D_{\zL'}$, $D_{\sT(M/G)}$, and $\dt D_{M/G}$ on
$\T\T(M/G)$. }
\end{ex}

\section{Acknowledgments}

The authors wish to thank F.~Przytycki for helpful discussions on dynamical systems.

%%%%%%%%%%%%%%%%%%%%%%%%%%%%%%%%%%%%%%%%%%%%%%%%
%%%%%%%%%%%%%%%%%%%%%%%%%%%%%%%%%%%%%%%%%%%%%%%%%%%%%%%%%%%%%o

\bigskip
\noindent Janusz Grabowski\\ Institute of Mathematics, Polish Academy of
Sciences\\\'Sniadeckich 8, P.O. Box 21, 00-956 Warszawa,
Poland\\{\tt jagrab@impan.gov.pl}\\\\
\noindent Miko\l aj Rotkiewicz\\
Institute of Mathematics,
                University of Warsaw \\
                Banacha 2, 02-097 Warszawa, Poland \\
                 {\tt mrotkiew@mimuw.edu.pl}
\end{document}